\documentclass[11pt]{article}

\usepackage{microtype}
\usepackage{amssymb}
\usepackage{graphicx}
\usepackage{subfigure}
\usepackage[usenames]{color}
\usepackage{colortbl}

\usepackage{enumitem}

\usepackage{latexsym}
\usepackage{fullpage, color}
\usepackage{array}
\usepackage{multirow}
\usepackage{tabularx}
\usepackage{wrapfig}
\usepackage{tablefootnote}
\usepackage[dvipsnames]{xcolor}
\usepackage{authblk}
\usepackage{algorithm,algorithmic}

\usepackage[utf8]{inputenc} 
\usepackage[T1]{fontenc}    
\usepackage{hyperref}       
\usepackage{url}          
\usepackage{booktabs}       
\usepackage{amsfonts}       
\usepackage{nicefrac}       
\usepackage{microtype} 
\usepackage{amsmath}
\usepackage{amsthm}
\usepackage{adjustbox}

\usepackage[capitalize,noabbrev]{cleveref}

\theoremstyle{plain}
\newtheorem{theorem}{Theorem}[section]

\newtheorem{lemma}[theorem]{Lemma}

\theoremstyle{definition}

\theoremstyle{remark}
\newtheorem{remark}[theorem]{Remark}

\usepackage[textsize=tiny]{todonotes}

\def\ag#1{{\color{black}#1}}
\def\an#1{{\color{black}#1}}
\def\dk#1{{\color{black}#1}}

\title{The Power of First-Order Smooth Optimization for Black-Box Non-Smooth Problems}

\author[1,2,3]{Alexander Gasnikov} \author[1,2]{Anton Novitskii} \author[1]{Vasilii Novitskii} 
\author[3]{Farshed Abdukhakimov} \author[3,1]{Dmitry Kamzolov} \author[1,4,3]{Aleksandr Beznosikov} 
\author[3]{Martin Tak\'a\v{c}}  \author[5]{Pavel Dvurechensky} \author[3]{Bin Gu}
\affil[1]{Moscow Institute of Physics and Technology, Dolgoprudny, Russia}
\affil[2]{ISP RAS Research Center for Trusted Artificial Intelligence, Moscow, Russia}
\affil[3]{Mohamed bin Zayed University of Artificial Intelligence, Abu Dhabi, UAE}
\affil[4]{National Research University Higher School of Economics, Moscow, Russian Federation}
\affil[5]{Weierstrass Institute for Applied Analysis and Stochastics, Berlin, Germany}

\begin{document}

\maketitle

\begin{abstract}

Gradient-free/zeroth-order methods for black-box convex optimization have been extensively studied in the last decade with the main focus on oracle \dk{call} complexity. In this paper, besides the oracle complexity, we focus also on iteration complexity, and propose a generic approach that, based on optimal first-order methods, allows to obtain in a black-box fashion new zeroth-order algorithms for non-smooth convex optimization problems. Our approach not only leads to optimal oracle complexity, but also allows to obtain iteration complexity similar to first-order methods, which, in turn, allows to exploit parallel computations to accelerate the convergence of our algorithms. We also elaborate on extensions for stochastic optimization problems, saddle-point problems, and distributed optimization.

\end{abstract}

\section{Problem Formulation}
We consider optimization problem
\begin{equation}\label{problem}
    \min_{x\in Q\subseteq \mathbb{R}^d} f(x)
\end{equation}
in the setting of a zeroth-order oracle. This means that an oracle returns the value $f(x)$ at a requested point $x$ \cite{conn2009introduction}, possibly with some adversarial noise that is uniformly bounded by a small value $\Delta >0$.
Let $\gamma > 0$ be a small number to be defined later and $Q_{\gamma}:=Q + B^d_2(\gamma)$, where $B^d_2(\gamma)$ is the Euclidean ball with center at $0$ and radius $\gamma > 0$ in $\mathbb{R}^d$.
Using these objects, we make the following assumptions.\footnote{Note that, for most of the algorithms in this paper, we can make these assumptions only on the intersection of $Q_{\gamma}$ and the ball $x^0+B^d_p(R)$ for some $p \in [1,2]$, where $x^0$ is the starting point of the algorithm and $R=O\left(\|x^0 - x_*\|_p\ln d\right)$ with $x_*$ being a solution of \eqref{problem} closest to $x^0$ \cite{gorbunov2019optimal}.}
\begin{itemize}[noitemsep,topsep=0pt,leftmargin=10pt]
    \item The set $Q$ is convex and the function $f$ is convex on the set $Q_{\gamma}$; \footnote{
    This assumption on the availability of  the objective values $f$ in a small vicinity $Q_{\gamma}$ of the feasible set $Q$  is quite standard in the literature, see, e.g., \cite{YOUSEFIAN201256,duchi2012randomized} and can be established in two ways.
    The first one is changing the set $Q$ in problem \eqref{problem} to a slightly smaller set $\widetilde{Q}$ such that $\widetilde{Q} + B^d_2(\gamma) \subseteq Q$, see, e.g., \cite{beznosikov2020gradient}. The second one is the extension of $f$ to the whole space $\mathbb{R}^d$ with preserving the convexity and Lipschitz continuity \cite{risteski2016algorithms}. More precisely, by changing the objective to $f_{new}(x):= f\left(\text{proj}_Q(x)\right) + \alpha\min_{y\in Q}\|x - y\|_2$.}
    \item The function $f$ is Lipschitz-continuous with constant $M$, i.e. $|f(y) - f(x) | \leq M \|y - x \|_p$ on $Q_{\gamma}$, where $p\in[1,2]$ and $\|\cdot\|_p$ is the $p$-norm. If $p=2$ we use the notation $M_2$ for the Lipschitz constant. 
 
\end{itemize}

This class of problems was widely investigated and optimal algorithms in terms of the number of zeroth-order oracle calls were developed in non-smooth setting \cite{nesterov2017random,duchi2015optimal,gasnikov2017stochastic,shamir2017optimal,bayandina2018gradient} and smooth setting \cite{nesterov2017random,dvurechensky2021accelerated,gorbunov2018accelerated}.  
At the same time, to the best of our knowledge, the development of optimal algorithms in terms of the number of iterations is still an open research question \cite{duchi2015optimal,bubeck2019complexity}.
The goal of this paper is to propose a generic approach that allows \dk{us} to construct algorithms with the best iteration complexity among the algorithms that have optimal oracle complexity. This can be seen as a two-criteria optimization problem if one would like to improve both: oracle complexity and iteration complexity. An important observation of this paper is that there is no need to sacrifice oracle complexity to obtain better iteration complexity.
A special focus is made on the possibility to use non-Euclidean geometry to define the algorithms.
\\[2pt]
The simplest and most illustrative result obtained in this work is as follows. Our generic approach allows to solve problem \eqref{problem} in the Euclidean geometry in $O(d^{1/4}M_2R/\varepsilon)$ successive iterations, each requiring $O(d^{3/4}M_2R/\varepsilon)$ oracle calls per iteration that can be made in parallel to make the total working time smaller. The dependence $\sim d^{1/4}/\varepsilon$ corresponds to the first part of the lower bound for iteration complexity \cite{diakonikolas2020lower,bubeck2019complexity}  $\min\left\{d^{1/4}\varepsilon^{-1}, d^{1/3}\varepsilon^{-2/3}\right\}$. Note that  this lower bound is obtained for a wider class of algorithms  that allow $\text{Poly}\left(d,\varepsilon^{-1}\right)$ oracle calls per iteration. On the contrary, our algorithm makes minimal possible number of oracle calls in each iteration in such a way that the total number of oracle calls is optimal.
\\[2pt]
Due to the page limitation, we next describe the main results and technical details are deferred to the Appendices.

\section{Smoothing Scheme}\label{SmSc}
In this section we describe our main scheme that allows \dk{us} to develop batch-parallel gradient-free methods for non-smooth convex problems based on batched-gradient algorithms for smooth stochastic convex problems. 
In the following sections, we generalize this scheme to non-smooth stochastic convex optimization problems and convex-concave saddle-point problems, including the problems with finite-sum structure.
\\[2pt]
The first element of our approach is the randomized \textit{smoothing} for non-smooth objective $f$. This approach is rather standard, and goes back to 1970s  \cite{ermoliev1976stochastic,nemirovsky1983problem,spall2005introduction}.
The smooth approximation to $f$ is defined as the function
$$
f_{\gamma}(x)=  \mathbb{E}_{u} f(x + u),
$$
where $u \sim RB^d_2(\gamma)$, i.e. $u$ is random vector uniformly distributed on $B^d_2(\gamma)$. The following theorem is a generalization of the results \cite{YOUSEFIAN201256,duchi2012randomized} for non-Euclidean norms.

\begin{theorem}[properties of $f_{\gamma}$]
\label{main_properties} For all $x,y\in Q$, we have
\begin{itemize}[noitemsep,topsep=0pt,leftmargin=10pt]
    \item $f(x) \leq f_{\gamma}(x) \leq f(x) + \gamma M_{\ag{2}}$;
    
    \item $f_{\gamma}(x)$ is $M$-Lipschitz: $$|f_{\gamma}(y) - f_{\gamma}(x) | \leq M \| y - x\|_p;$$
    
    \item $f_{\gamma}(x)$ has $L = 
    \dfrac{\sqrt{d} M}{\gamma} $-Lipschitz gradient:
    $$\|\nabla f_{\gamma}(y) - \nabla f_{\gamma}(x) \|_q \leq L \| y - x\|_p,$$ 
     where $q$ us such that $1/p + 1/q = 1$. 
\end{itemize}
\end{theorem}
See Appendix \ref{main_properties_proof} for the proof. 
\\[1pt]
The second very important element of our approach goes back to \cite{shamir2017optimal}, who proposes a special unbiased stochastic gradient for $f_{\gamma}$ with small variance:
\begin{equation}\label{sg}
    \nabla f_{\gamma}(x,e) = d\frac{f(x+\gamma e) - f(x-\gamma e)}{2\gamma}e,
\end{equation}
where $e \sim RS^d_2(1)$ is a random vector uniformly distributed on $S^d_2(1)$ -- Euclidean unit sphere with center at $0$ in $\mathbb{R}^d$. To simplify the further derivations, here and below we make a slight abuse of notation and denote by $f$ the value returned by the oracle, which can be an approximation to the value of the objective up to a small error, bounded by  $\Delta$.

We note that an alternative way to define a stochastic approximation to $\nabla f_{\gamma}(x)$ is based on the double smoothing technique of B. Polyak \cite{polyak1987introduction,bayandina2018gradient}. This approach is more complicated and requires stronger assumptions on the noise $\Delta$ (see Theorem~\ref{main_properties2}). Thus, we use the approach of \cite{shamir2017optimal}.

The following theorem is a combination of the results from \cite{shamir2017optimal,gorbunov2019upper,beznosikov2020derivative}.
    \begin{theorem}[properties of $\nabla f_{\gamma}(x,e)$]
\label{main_properties2} 
For all $x\in Q$, we have
 \begin{itemize}[noitemsep,topsep=0pt,leftmargin=10pt] 
    \item $\nabla f_{\gamma}(x,e)$ is an unbiased approximation for $\nabla f_{\gamma}(x)$:\footnote{For simplicity, we assume here (and everywhere where we talk about unbiased estimates) that the small noise in the function value is random (but not not necessary i.i.d.) with zero mean. For the general setting see \cite{beznosikov2020derivative} \ag{and Appendix~\ref{noisy}}.}
    $\mathbb{E}_e \left[\nabla f_{\gamma}(x, e) \right] = \nabla f_{\gamma}(x)$;
    
    \item  $\nabla f_{\gamma}(x,e)$ has bounded variance (second moment):
    $$ \mathbb{E}_e \left[ \| \nabla f_{\gamma}(x, e) \|^2_q \right] \leq 
    \kappa(p,d)\cdot\left(d M_2^2 + \dfrac{d^2\Delta^2}{\gamma^2} \right),$$
   where $1/p + 1/q = 1$ and
   \begin{equation*}
   \begin{aligned}
\kappa(p,d) &= O\left(\sqrt{\mathbb{E}_e\|e\|_q^4}\right) = \ag{O\left(\min\left\{q,\ln d\right\}d^{2/q - 1}\right)}\\
&=
 \begin{cases}
   O(1), \;  p = 2~\ag{(q=2)}\\
   O\left((\ln d)/d\right), \; p = 1~\ag{(q=\infty)}.
 \end{cases}
   \end{aligned}
\end{equation*}
Moreover, if $\Delta$ is sufficiently small, then 
 $$\mathbb{E}_e \left[ \| \nabla f_{\gamma}(x, e) \|^2_q \right]\leq 2\kappa(p,d)dM_2^2.$$
\end{itemize}

\end{theorem}

See Appendix \ref{theorem_22_proof} for the proof.
\begin{remark}\label{stoch}
If the zeroth-order oracle returns an unbiased noisy stochastic  function value $f(x,\xi)$ ($\mathbb{E}_{\xi} f(x,\xi) = f(x)$), then with two-point oracle we can introduce the following counterpart of \eqref{sg}
\begin{equation}
\label{sg_stoch}
     \nabla f_{\gamma}(x,\xi,e) = d\frac{f(x+\gamma e,\xi) - f(x-\gamma e,\xi)}{2\gamma}e.
\end{equation}
 Theorem~\ref{main_properties2} remains valid with the appropriate changes of $\nabla f_{\gamma}(x,e)$  to $\nabla f_{\gamma}(x,\xi,e)$, the expectation $\mathbb{E}_{e}$ to the expectation $\mathbb{E}_{e,\xi}$, and the redefinition of $M_2$ as a constant satisfying $\mathbb{E}_{\xi}\|\nabla_x f(x,\xi)\|_2^2\le M^2_2$ for all $x\in Q_{\gamma}$.
\end{remark}
Based on the two elements above, we are now in a position to describe our general approach, which for shortness, we refer to as the \textit{Smoothing scheme} (technique).
\\[2pt]
Assume that we have some batched algorithm \textbf{A}($L,\sigma^2$) that solves problem \eqref{problem} under the assumption that $f$ is smooth and satisfies
\begin{equation}\label{LL}
\|\nabla f(y) - \nabla f(x) \|_q \leq L \| y - x\|_p, \; \forall x,y\in Q_{\gamma},
\end{equation}
and by using \dk{a} stochastic first-order oracle that depends on a random variable $\eta$ and returns at a point $x$ an unbiased stochastic gradient $\nabla_x f(x,\eta)$ with bounded variance:  
\begin{equation}
\label{sigma^2}
     \mathbb{E}_{\eta} \left[ \| \nabla_x f(x, \eta) - \nabla f(x)\|^2_q \right]\le \sigma^2.
\end{equation}
Further, we assume that, to reach $\varepsilon$-suboptimality in expectation, this algorithm requires $N(L,\varepsilon)$ successive iterations and $T(L,\sigma^2,\varepsilon)$ stochastic first-order oracle calls, i.e. \textbf{A}($L,\sigma^2$) allows batch-parallelization with the average batch size $B(L,\sigma^2,\varepsilon) = T(L,\sigma^2,\varepsilon)/N(L,\varepsilon)$.
 \\[2pt]
Our approach consists of applying \textbf{A}($L,\sigma^2$) to the smoothed problem
 \begin{equation}\label{sm_problem}
    \min_{x\in Q\subseteq \mathbb{R}^d} f_{\gamma}(x)
\end{equation}
with 
\begin{equation}\label{gamma}
  \gamma = \varepsilon/(2M_{\ag{2}})  
\end{equation}
and $\eta = e$, $\nabla_x f(x, \eta) = \nabla f_{\gamma}(x, e)$, where $\varepsilon > 0$ is the desired accuracy for solving problem \eqref{problem} in terms of the suboptimality expectation.
\\[2pt]
According to Theorem~\ref{main_properties}, an $(\varepsilon/2)$-solution to \eqref{sm_problem} is an $\varepsilon$-solution to the initial problem \eqref{problem}.
According to Theorem~\ref{main_properties} and \eqref{gamma} we have 
\begin{equation}\label{L}
   L \le \frac{2\sqrt{d}M M_{\ag{2}}}{\varepsilon}, 
\end{equation}
and, according to Theorem~\ref{main_properties2}, we have 
\begin{equation}\label{var}
  \sigma^2 \le 2\kappa(p,d)dM_2^2  
\end{equation}
if $\Delta$ is sufficiently small.
\\[2pt]
Thus, we obtain that \textbf{A}($L,\sigma^2$) implemented using stochastic gradient \eqref{sg} is a zeroth-order method for solving non-smooth problem \eqref{problem}. Moreover, to solve problem \eqref{problem} with accuracy $\varepsilon>0$ this method suffices to make 
$$
N\left(\frac{2\sqrt{d}M M_{\ag{2}}}{\varepsilon},\varepsilon\right) \; \text{ successive iterations and }
$$
$$
2T\left(\frac{2\sqrt{d}M M_{\ag{2}}}{\varepsilon},2\kappa(p,d)dM_2^2,\varepsilon\right) \; \text{ zeroth-order oracle calls. }
$$

We underline that this approach is flexible and generic as we can take different algorithms as \textbf{A}($L,\sigma^2$). For example, if we take batched Accelerated gradient method \ag{\cite{cotter2011better,lan2012optimal,devolder2013exactness,dvurechensky2016stochastic,gorbunov2019optimal}}, then from \eqref{L}, \eqref{var} we have that
\begin{center}
$N(L,\varepsilon) = O\left(\sqrt{\frac{LR^2}{\varepsilon}}\right)=O\left(\frac{d^{1/4}\ag{\sqrt{M M_2}}R}{\varepsilon}\right),$ 
$T(L,\sigma^2,\varepsilon) = \tilde{O}\left(\max\left\{N(L,\varepsilon),\frac{\sigma^2R^2}{\varepsilon^2}\right\}\right)$
$=\tilde{O}\left(\frac{\kappa(p,d)d M_2^2R^2}{\varepsilon^2}\right),$ 
\end{center}
\dk{where $\tilde{O}\left( \, \cdot \, \right)$ is a convergence rate up to a logarithmic factor}
Here $R = O(\|x^0 - x_*\|_p\ln d)$ with $x^0$ being the starting point and  $x_*$ being the solution to \eqref{problem} closest to $x^0$. The last equality assumes also that $\varepsilon \lesssim d^{-1/4}M_2^{\ag{3/2}}R/M^{\ag{1/2}}$ when $p=1$.

\begin{theorem}\label{FGM}
Based on the batched Accelerated gradient method, the Smoothing scheme applied to non-smooth problem \eqref{problem}, provides a gradient-free method with $$
O\left(\frac{d^{1/4}\ag{\sqrt{M M_2}}R}{\varepsilon}\right) \; \text{ successive iterations and}
$$ 
$$\tilde{O}\left(\frac{\kappa(p,d)d M_2^2R^2}{\varepsilon^2}\right)=
\begin{cases}
   \tilde{O}\left(\frac{d M_2^2R^2}{\varepsilon^2}\right), \;  p = 2\\
   \tilde{O}\left(\frac{(\ln d) M_2^2R^2}{\varepsilon^2}\right), \; p = 1.
 \end{cases}
 $$
zeroth-order oracle calls.
\end{theorem}

See Appendix \ref{th2.4proof} for the proof.

Next we make several remarks on some related works. First, an important difference between our \textit{Smoothing scheme} and the work of \cite{duchi2015optimal} is that, unlike them we do not assume the smoothness of the objective $f$ and we use a different finite-difference approximation \eqref{sg} due to \cite{shamir2017optimal}. 
Second, in \cite{scaman2019optimal}, the authors use a similar smoothing technique to reduce the number of communications in distributed non-smooth convex optimization algorithms. Unlike their exact first-order oracle setting, we consider zeroth-order oracle model and construct stochastic approximation to the gradient of the smoothed function.
Finally, the authors of \cite{bubeck2019complexity} propose a close technique with an accelerated higher-order method \cite{nesterov2021implementable,gasnikov2019near,agafonov2020inexact} playing the role of \textbf{A}($L,\sigma^2$). For the particular setting of $p=2$, they obtain a better in some regimes bound $N\sim d^{1/3}/\varepsilon^{2/3}$ for the number of iterations. Yet, they have significantly worse oracle complexity $T$, which makes it unclear how to use their results in practice.   
To sum up, despite some similarities with other works, the proposed \textit{Smoothing scheme} is, to the best of our knowledge, new,  general, and flexible. Moreover, as we show below it is quite universal and can be applied to many different problems.

\section{Applications Of the Smoothing Scheme}
\subsection{Stochastic Optimization}\label{SO}
Based on Remark~\ref{stoch}, we can consider the non-smooth stochastic convex optimization problem:
\begin{equation}\label{stoch_problem}
    \min_{x\in Q\subseteq \mathbb{R}^d} \left\{f(x):=\mathbb{E}_{\xi} f(x,\xi)\right\}
\end{equation}
with \dk{a} two-point zeroth-order oracle that returns the values $\{f(x_i,\xi)\}_{i=1}^2$ given two points $x_1,x_2$.
The result of Theorem~\ref{FGM} still holds\footnote{Note, that for the stochastic optimization problem \eqref{stoch_problem} it is important that \textbf{A}($L,\sigma^2$) requires $L$ to be defined according to \eqref{LL}, rather than as $L$ satisfying $$\|\nabla f(y,\xi) - \nabla f(x,\xi) \|_q \leq L \| y - x\|_p,$$ for all
$x,y\in Q_{\gamma}$
and all $\xi$. For example, this means that we can not apply the \textit{Smoothing technique} to batched Accelerated gradient method with interpolation \cite{woodworth2021even}.} 
if $M_2$ is redefined to be a constant satisfying $\mathbb{E}_{\xi}\|\nabla_x f(x,\xi)\|_2^2\le M^2_2$ for all
$x\in Q_{\gamma}$, $\eta$ is set to be the pair $(\xi,e)$. \dk{$e$ and $\xi$ are independent between iterations with available samples.}
Moreover, by using the stochastic gradient clipping \cite{gorbunov2020stochastic}, we can prove a stronger result and guarantee $\varepsilon$-suboptimality with high-probability (with exponential concentration) independently of distribution of $\nabla_x f(x,\xi)$.

If the two-point feedback as in \eqref{sg_stoch} is not available, our \textit{Smoothing technique} can utilize the one-point feedback by using the unbiased estimate \cite{nemirovsky1983problem,flaxman2005online,gasnikov2017stochastic}:
\begin{equation*}
    \nabla f_{\gamma}(x,\xi,e) = d\frac{f(x+\gamma e,\xi)}{\gamma}e,
\end{equation*}
with \cite{gasnikov2017stochastic}
$$ \mathbb{E}_{\xi,e} \left[ \| \nabla f_{\gamma}(x, \xi, e) \|^2_q \right] \le   \begin{cases}
   \frac{(q-1)d^{1+2/q} G^2}{\gamma^2}, \;  q\in[2,2\ln d]\\
  \frac{4d(\ln d) G^2}{\gamma^2}, \; q \in (2\ln d,\infty),
 \end{cases}$$
where $\gamma$ is defined in \eqref{gamma} and it is assumed that  $\mathbb{E}_{\xi} \left[ | f(x, \xi) |^2\right] \le G^2$ for all 
$x\in Q_{\gamma}$.
Thus, the \textit{Smoothing technique} can be generalized to the one-point feedback setup by replacing the RHS of \eqref{var} by the above estimate. This leads to the same iteration complexity $N$, but increases the oracle complexity $T$, and, consequently, the batch size  $B$ at each iteration.  

\subsection{Finite-sum Problems}\label{STP}
As a special case of \eqref{stoch_problem} with $\xi$ uniformly distributed on $1,...,m$, we can consider the finite-sum (Empirical Risk Minimization) problem 
\begin{equation}\label{sum}
    \min_{x\in Q\subseteq \mathbb{R}^d} f(x):=\mathbb{E}_{\xi} f(x,\xi) = \frac{1}{m}\sum_{k=1}^m f_k(x).
\end{equation}
Clearly, if we have incremental zeroth-order oracle, i.e. zeroth-order oracle for each $f_k$, we are in the setting of two-point feedback in the sense of Section~\ref{SO}. Thus, we can apply Theorem~\ref{FGM}, where $M_2$ is defined as $\max_{k=1,...,m}\|\nabla f_k(x)\|_2\le M_2$ for all 
$x\in Q_{\gamma}$. 

At the same time, problem \eqref{sum} has specific structure, that allows to split batching among nodes (to make algorithm centralized distributed among $m$ nodes) if batch size $$B = \frac{T}{N} \simeq \frac{dM_2^2R^2/\varepsilon^2}{d^{1/4}M_2R/\varepsilon}= d^{3/4}\frac{M_2R}{\varepsilon}$$ is greater than $m$, i.e.\footnote{For clarity in this subsection we consider the Euclidean case with $p=2$.} $m \lesssim d^{3/4}M_2R/\varepsilon$. As it is known, in Machine Learning applications $m$ can be as large as $O\left(dM_2^2R^2/\varepsilon^2\right)$ \cite{shapiro2005complexity,feldman2016generalization} or $O\left(M_2^2R^2/\varepsilon^2\right)$, if a proper regularization is applied \cite{shalev2009stochastic,shalev2014understanding}. In both cases $m$ can be very large. 

So, in this case ($d^{3/4}M_2R/\varepsilon \lesssim m \lesssim d^{3/2}M_2^2R^2/\varepsilon^2$), to preserve the total oracle complexity one can use stochastic Accelerated Variance Reduced algorithms \cite{lan2018random,kulunchakov2019estimate,kulunchakov2020optimisation,lan2020first} with (see \eqref{L}, \eqref{var})
\begin{center}
$N(L,\varepsilon) = \tilde{O}\left(m + \sqrt{\frac{mLR^2}{\varepsilon}}\right)=O\left(\frac{d^{1/4}\sqrt{m}M_2R}{\varepsilon}\right),$ $T(L,\sigma^2,\varepsilon) = \tilde{O}\left(\max\left\{N(L,\varepsilon),\frac{\sigma^2R^2}{\varepsilon^2}\right\}\right)=$
$=\tilde{O}\left(\frac{dM_2^2R^2}{\varepsilon^2}\right).$ 
\end{center}
The number of successive iterations grows, but now $m$-nodes distribution of batching is possible if $$m\lesssim B = \frac{T}{N} \simeq \frac{dM_2^2R^2/\varepsilon^2}{d^{1/4}\sqrt{m}M_2R/\varepsilon}=d^{3/4}\frac{M_2R}{\sqrt{m}\varepsilon},$$
that is
$$m\lesssim d^{3/2}\frac{M_2^2R^2}{\varepsilon^2}.$$
This regime is quite natural for Machine Learning and Statistical applications \cite{shalev2014understanding,shapiro2021lectures}.
\\[2pt]
The  results of this subsection can be generalized to the stochastic optimization problems with $f_k(x):=\mathbb{E}_{\xi} f_k(x,\xi)$, see (Section~\ref{SO}).

\subsection{Strongly Convex Problems}\label{SC}
By using the restart technique \cite{nemirovsky1983problem,juditsky2014deterministic} we can prove a counterpart of Theorem~\ref{FGM} for the case when $f$ is $\mu$-strongly convex w.r.t. the $p$-norm for some $p\in[1,2]$ and $\mu \ge \varepsilon/R^2$.
\begin{theorem}\label{SFGM}
Based on the batched Accelerated gradient method, the Smoothing scheme applied to non-smooth and strongly convex problem \eqref{problem}, provides a gradient-free method with
$\tilde{O}\left(\frac{d^{1/4}\ag{\sqrt{M M_2}}}{\sqrt{\mu\varepsilon}}\right)$   successive iterations and   $\tilde{O}\left(\frac{\kappa(p,d)d M_2^2}{\mu\varepsilon}\right)$
zeroth-order oracle calls, where $\kappa(p,d)$ \dk{is bounded as} in Theorem~\ref{FGM}. 
\\[2pt]
Moreover, the same holds for stochastic optimization problem \eqref{stoch_problem} if $M_2$ is defined as $\mathbb{E}_{\xi}\|\nabla_x f(x,\xi)\|_2^2\le M^2_2$ for all
$x\in Q_{\gamma}$.
\end{theorem}
See Appendix \ref{th3.1.proof} for the proof. 

\subsection{Saddle-point Problems}\label{S-PP}
In this subsection we consider non-smooth convex-concave saddle-point problem
\begin{equation}\label{SPP}
    \min_{x\in Q_x\subseteq \mathbb{R}^{d_x}} \max_{y\in Q_y\subseteq \mathbb{R}^{d_y}} f(x,y).
\end{equation}
Gradient-free methods for convex-concave saddle-point problems were studied in \cite{beznosikov2020gradient,beznosikov2021one,gladin2021solving,sadiev2021zeroth} with the main focus on the complexity in terms of the number of zeroth-order oracle calls. Unlike these papers, we focus here also on the iteration complexity. 

Applying \textit{Smoothing technique} separately to $x$-variables and $y$-variables, we obtain almost the same results as for optimization problems with the only difference in Theorem~\ref{main_properties}: instead of
$$f(x) \leq f_{\gamma}(x) \leq f(x) + \gamma M_{\ag{2}}$$ we have $$f(x\ag{, y}) - \gamma_y M_{\ag{2,}y}\leq f_{\gamma}(x\ag{, y}) \leq f(x\ag{, y}) + \gamma_x M_{\ag{2,}x}.$$
This leads to a clear counterpart of \eqref{gamma} for choosing $\gamma = \left(\gamma_x,\gamma_y\right)$, where $M_{\ag{2,}x}$, $M_{\ag{2,}y}$ -- corresponding Lipschitz constants in \ag{$2$-norm}.

If we take as \textbf{A}($L,\sigma^2$) the batched Mirror-Prox or the batched Operator extrapolation method or the batched Extragradient method \cite{juditsky2011solving,kotsalis2020simple,gorbunov2021stochastic}, using  \eqref{L}, \eqref{var}, we obtain the following bounds
\begin{center}
$N(L,\varepsilon) = O\left(\frac{LR^2}{\varepsilon}\right)=O\left(\frac{\sqrt{d}\ag{M M_2}R^2}{\varepsilon^2}\right),$ $ T(L,\sigma^2,\varepsilon) = O\left(\max\left\{N(L,\varepsilon),\frac{\sigma^2R^2}{\varepsilon^2}\right\}\right)$
$=O\left(\frac{\kappa(p,d)d M_2^2R^2}{\varepsilon^2}\right),$ 
\end{center}
where $d = \max\{d_x,d_y\}$,
\ag{$M_2=\max\left\{M_{2,x},M_{2,y}\right\}$,} $R$ depends on the criteria. For example, if $\varepsilon$ is expected accuracy in fair duality gap \cite{juditsky2011solving}, then $R$ is a diameter in $p$-norm of $Q_x\otimes Q_y$ up to a $\ln d$-factor\dk{, where $\otimes$ is a Cartesian product of two sets}. The last equality assumes that $d \lesssim (M_2/M)^{\ag{2}}$ when $p=1$. This result is also correct for stochastic saddle-point problems with proper redefinition of what is $M_2$, see Section~\ref{SO}.
\\[2pt]
We see that due to the lack of acceleration for saddle-point problems the batch-parallelization effect is much more modest than for convex optimization problems.
\\[2pt]
By using the restart technique we can generalize these results to $\mu$-strongly convex, $\mu$-strongly concave case, see Section~\ref{SC}. Alternatively, we can combine the \textit{Smoothing technique} with Stochastic Accelerated Primal-Dual method from \cite{zhang2021robust} for \eqref{SPP} with $f(x,y)$ being $\mu_x$-strongly convex and $\mu_y$-strongly concave in $2$-norm (Euclidean setup). In this case we obtain the following bounds
\begin{center}
$N(\{L\},\varepsilon) = \tilde{O}\left(\frac{L_{xx}}{\mu_x} + \frac{\max\left\{L_{xy},L_{yx}\right\}}{\sqrt{\mu_x\mu_y}}+ \frac{L_{yy}}{\mu_y} \right)$
$=\tilde{O}\left(\frac{\sqrt{d_x}M_{2,x}^2}{\mu_x\varepsilon}+ \frac{M_{2,x}M_{2,y}\max\left\{\sqrt{d_x},\sqrt{d_y}\right\}}{\sqrt{\mu_x\mu_y}\varepsilon}  + \frac{\sqrt{d_y}M_{2,y}^2}{\mu_y\varepsilon}\right),$
$T(\{L\},\{\sigma^2\},\varepsilon) = \tilde{O}\left(\max\left\{N(\{L\},\varepsilon),\frac{\sigma_x^2}{\mu_x\varepsilon} + \frac{\sigma_y^2}{\mu_y\varepsilon}\right\}\right)$
$=\tilde{O}\left(\frac{d_x M_{2,x}^2}{\mu_x\varepsilon}+\frac{d_y M_{2,y}^2}{\mu_y\varepsilon}\right),$ 
\end{center}
where we use subscripts corresponding to $x$ or $y$ variables, e.g. \dk{$L_{xx}, L_{xy}, L_{yx}, L_{yy}$ is defined as
$$\|\nabla_x f_{\gamma}(x_2,y) - \nabla_x f_{\gamma}(x_1,y)\|_2 \le L_{xy}\|x_2 - x_1\|_2$$
$$\|\nabla_x f_{\gamma}(x,y_2) - \nabla_x f_{\gamma}(x,y_1)\|_2 \le L_{xy}\|y_2 - y_1\|_2$$
$$\|\nabla_y f_{\gamma}(x_2,y) - \nabla_y f_{\gamma}(x_1,y)\|_2 \le L_{yx}\|x_2 - x_1\|_2$$
$$\|\nabla_y f_{\gamma}(x,y_2) - \nabla_y f_{\gamma}(x,y_1)\|_2 \le L_{yy}\|y_2 - y_1\|_2$$}
for all $(x,y)\in Q_{x,\gamma_x}\otimes Q_{y,\gamma_y}$.
The only non-trivial calculation here is estimation of $L_{xy}$, $L_{yx}$, see Appendix \ref{L_xy_proof} for the proof. 
The other constants $\{L\}$ and $\{\sigma^2\}$ are defined according to the standard \textit{Smoothing scheme} with variables $x$ or $y$  corresponding to subscripts.
\\[2pt]
Note, that most of the results for saddle-point problems (i.e. mentioned result from \cite{zhang2021robust} or finite-sum composite generalization \cite{vladislav2021accelerated}) with different constants of smoothness and strong convexity/concavity were obtained based on \dk{the} Accelerated gradient method for convex problems and Catalyst envelope, that allows \dk{us} to generalize it to saddle-point problems \cite{pmlr-v125-lin20a}. There exist also loop-less (direct) accelerated methods that save $\ln(\varepsilon^{-1})$-factor in the complexity, for $\mu_x$-strongly convex, $\mu_y$-strongly concave saddle-point problems \cite{kovalev2021accelerated}. But even for composite bilinear saddle-point problems (with different smoothness and strong convexity constants) there is still a gap between the state-of-the-art upper bounds \cite{kovalev2021accelerated} and lower bounds \cite{zhang2021lower}.
 
\subsection{Distributed Optimization}\label{DO}
In decentralized distributed convex optimization and convex-concave saddle-point problems optimal methods (both in terms of communication rounds and oracle calls) were developed in the Euclidean setup, see, e.g.  surveys \cite{gorbunov2020recent,Dvinskikh2021}. In particular, there exists  a batched-consensus-projected Accelerated gradient method \cite{rogozin2021accelerated} that, for $\mu$-strongly convex in $2$-norm $f$ from \eqref{sum} with $L$-Lipschitz gradient in $2$-norm, requires  
$$N(L,\varepsilon) = \tilde{O}\left(\sqrt{\chi\frac{L}{\mu}}\right)=\tilde{O}\left(\frac{\sqrt{\chi}d^{1/4}M_2}{\sqrt{\mu\varepsilon}}\right)$$
communication rounds and
$$T(L,\sigma^2,\varepsilon) = \tilde{O}\left(\max\left\{N(L,\varepsilon),\frac{\sigma^2}{\mu\varepsilon^2}\right\}\right)=\tilde{O}\left(\frac{d M_2^2}{\mu\varepsilon}\right)$$
oracle calls per node, where $M_2$ is defined in Section~\ref{STP}, $\chi$ -- condition number of the Laplace matrix of communication network or square of worst-case condition number for time-varying networks \cite{rogozin2021accelerated}. 
\\[2pt]
This batched-consensus-projected Accelerated gradient method is optimal (up to a $\ln (\varepsilon^{-1})$-factor) as a decentralized method, but the \textit{Smoothing technique} provides a gradient-free method that is not the best (state-of-the-art) method in terms of communication rounds. The best one requires $\tilde{O}\left(\sqrt{\chi}M_2/\sqrt{\mu\varepsilon}\right)$ communication rounds \cite{beznosikov2020derivative}. This holds for stochastic decentralized convex problems with two-point feedback and one-point feedback. For the one-point feedback, optimal decentralized method is described in \cite{stepanov2021onepoint}. This method is better also by a $\sim d^{1/4}$-factor in terms of the number of communication rounds.
\\[2pt]
For saddle-point problems, by replacing in the \textit{Smoothing scheme}  the batched-consensus Accelerated gradient method \cite{rogozin2021accelerated}, which is optimal for decentralized convex problems, with the batched-consensus Extragradient method \cite{aleks2020distributed}, which is optimal for decentralized convex-concave saddle-point problems, we \dk{lose a} $\sim \sqrt{d}$-factor in the number of communication rounds in comparison with optimal gradient-free methods for non-smooth decentralized saddle-point problems.
\\[2pt]
To sum up, in distributed optimization, for the first time, we have a situation where the \textit{Smoothing scheme} generates a non-optimal method from an optimal one. 

\section{Discussion}\label{Discussion}
\subsection{Superposition Of Different Techniques}\label{Techniques}
In convex optimization and convex-concave saddle-point problems there are several generic techniques that allow black-box reduction of known methods to develop new methods for new problem classes \cite{gasnikov2017universal}.
For example, the \textit{Restart technique} mentioned in Section~\ref{SC} allows to construct methods for strongly convex problems based on methods for convex problems;
the \textit{Catalyst envelope} mentioned in Section~\ref{S-PP} allows to construct methods for convex-concave saddle-point problems based on methods developed for convex optimization problems; 
\textit{Batching technique} mentioned in Section~\ref{SmSc} allows to construct methods for stochastic problems based on methods developed for deterministic problems;
\textit{Consensus-projection technique} mentioned in Section~\ref{DO} allows to build decentralized distributed methods for convex problems based on non-decentralized methods developed for non-decentralized convex problems. \\[2pt]
The first important property of all these reduction techniques, informally speaking, is that all of them preserve the optimality of the method: an optimal method after applying any of the techniques becomes optimal \cite{gasnikov2017universal} for the new class of problems. We remark that the \textit{Catalyst envelope} leads to optimal algorithms only on a certain (yet large enough) class of problems, see Section~\ref{S-PP} for details. Also, when considering saddle-point problems with different strong convexity constants, it is better to use a direct method rather than the \textit{Restart technique}. 
Despite these limitations, the optimality-preserving property is very useful in practice.
\\[2pt]
The second important property is that a superposition or different combination of these reduction techniques also preserves the optimality of algorithms. This can be demonstrated via the batched-consensus-projected Accelerated gradient method from Section~\ref{DO}. This algorithm is obtained via a combination of the \textit{Batching technique} and  the \textit{Consensus-projection technique} applied to Accelerated gradient method. Moreover, the \textit{Restart technique} can be added to this combination in order to obtain the algorithm for strongly convex problems.
\\[2pt]
Note that in some cases these techniques require a proper generalization before they can be applied as a part of combination.
For example, when we generalize the batched-consensus-projected Accelerated gradient method to saddle-point problems, instead of the standard \textit{Catalyst envelope} developed in this context in \cite{pmlr-v125-lin20a}, we should use a special decentralized stochastic (batched) version of the \textit{Catalyst envelope}, that can be developed from \cite{kulunchakov2020optimisation,tian2021acceleration}.
\\[2pt]
We expect that the \textit{Smoothing scheme (technique)} developed in this paper will take its rightful place in the mentioned above (not exhaustive) list of useful reduction techniques that allow \dk{us} to develop new methods based on existing ones.
\\[2pt]
In the previous (sub)sections we demonstrated that the \textit{Smoothing scheme} is quite generic and allows a black box reduction of algorithms for smooth problems to solve non-smooth problems. Moreover, except the combination with \textit{Consensus-projection technique}, the \textit{Smoothing scheme} allows to obtain optimal algorithms for non-smooth black-box problems with zeroth-order oracle. \textbf{Thus, as one of our main contributions in this paper, we consider the observation that the \textit{Smoothing scheme} can be developed in a such a way, that it can be used in different combinations with other reduction techniques.}

\subsection{Batching Technique}\label{BT}
In the \textit{Smoothing scheme} an input algorithm should be a batched-gradient algorithm, and, thus, the \textit{Smoothing scheme}  strongly depends on the \textit{Batching technique}, which we now describe in more detail. In Section~\ref{SmSc} we mentioned some particular algorithms \cite{cotter2011better,lan2012optimal,devolder2013exactness,dvurechensky2016stochastic} which do not constitute a generic technique that allows \dk{us} to utilize  an arbitrary algorithm, which solves \dk{a} deterministic smooth convex problem.
Some attempts to propose a generic \textit{Batching technique} were made \cite{dvinskikh2020accelerated,gasnikov2017universal} in a much more general setting of inexact models of the objective, which we do \dk{not} consider here. Instead, we describe here a simple version in the convex case and for the Euclidean setting with $p=2$.
\\[2pt]
First of all, following \cite{devolder2013exactness,dvinskikh2020accelerated,Dvinskikh2021} we introduce the notion of $(\delta_1,\delta_2,L)$-oracle. We say that for the problem \eqref{problem} we have an access to $(\delta_1,\delta_2,L)$-oracle at a point $x$ if we can evaluate a vector $\nabla_{\delta} f(x)$ such that, for all $x,y \in Q_{\gamma}$,

$$ - \delta_1  \le f(y) - f(x) - \langle\nabla_{\delta} f(x), y-x \rangle  \le  \frac{L}{2}\|y-x\|_2^2 + \delta_2,$$
where $\mathbb{E}\delta_1 = 0$ ($\delta_1$ is independently taken at each oracle call),  $\mathbb{E}\delta_2 \le \delta$. 
Note that the left inequality corresponds to the definition of $\delta_1$-(sub)gradient \cite{polyak1987introduction} and reduces to the convexity property in the case $\delta_1 = 0$. In this case the LHS holds with $\nabla_{\delta} f(x) = \nabla f(x)$. The right inequality in the case when $\delta_2 = 0$ is a consequence\footnote{Note, that the right inequality in the case when $\delta_2 = 0$ is not equivalent to \eqref{LL}, but is typically sufficient to obtain optimal (up to constant factors) bounds on the rate of convergence  of different methods \cite{cite-key}.} of \eqref{LL}. Let us consider an algorithm \textbf{A}($L,\delta_1,\delta_2$) 
that converges with the rate\footnote{$N$ is a number of iterations which up to a constant factor is equal to the number of $ (\delta_1,\delta_2,L)$-oracle calls \ag{and $x_*$ is a solution of \eqref{problem}}. We can consider more specific rates of convergence for problems with additional structure and develop \textit{Batching technique} in a similar way.}
\begin{equation}\label{RC}
    \mathbb{E} f(x^N) - f(x_*) = O\left( \frac{LR^2}{N^{\alpha}} + N^{\beta}\delta \right).
\end{equation}
The \textit{batching technique}, applied to the problem \eqref{stoch_problem} with $L$-Lipschitz gradient (in 2-norm), is based on the use of the mini-batch stochastic approximation of the gradient
\[
\nabla_{\delta} f(x) = \frac{1}{r}\sum_{j=1}^r \nabla_x f(x,\xi^j)\\
\]
in \textbf{A}($L,\delta_1,\delta_2$), where $\{\xi^j\}_{j=1}^r$ are sampled independently and  $r$ is an appropriate batch size. The choice of $r$ is based on the following relations
\[\langle \nabla_{\delta} f(x) - \nabla f(x),y-x\rangle \le
\]
\[
\le \frac{1}{2L}\|\nabla_{\delta} f(x) - \nabla f(x)\|_2^2 + \frac{L}{2}\|y-x\|_2^2,
\]
\[
\mathbb{E}_{\{\xi^j\}}\left[ \| \nabla_{\delta} f(x) - \nabla f(x)\|_2^2 \right] \le \frac{\sigma^2}{r}, \]
where $\sigma^2$ is the variance of $\nabla_x f(x,\xi)$, see \eqref{sigma^2}. Hence, if
\[\delta \le \frac{1}{2L}
\max_{x\in Q_{\gamma}}
\mathbb{E}_{\{\xi^j\}_{j=1}^r}\left[ \| \nabla_{\delta} f(x) - \nabla f(x)\|_2^2 \right], \; \text{i.e.}\]
$\delta = \frac{\sigma^2}{2Lr}$, we have that  \textbf{A}($2L,\delta_1,\delta_2$) 
converges with the rate given in \eqref{RC}.
From \eqref{RC} we see that to obtain
 \[\mathbb{E} f(x^N) - f(x_*) \le \varepsilon \]
it suffices to take
\[N = O\left(\left(\frac{LR^2}{\varepsilon}\right)^{1/\alpha}\right) \quad \text{and}
\quad
r = O\left(\frac{\sigma^2N^{\beta}}{L\varepsilon}\right). 
\]
In particular, for the Accelerated gradient method we have that $\alpha = 2$, $\beta = 1$ \cite{devolder2013exactness,dvinskikh2020accelerated}. In this case, we obtain the complexity bounds for batched Accelerated gradient methods mentioned in Section~\ref{SmSc}:
\begin{center}
$N=O\left(\sqrt{LR^2/\varepsilon}\right)$, $B = r = O\left(\sigma^2R/\left(\sqrt{L}\varepsilon^{3/2}\right)\right)$, $T = N\cdot B = O\left(\sigma^2R^2/\varepsilon^2\right)$.
\end{center}
Note that, based on \cite{stonyakin2021}, the described above \textit{Batching technique} can be applied to saddle-point problems. In particular, this allows to obtain new methods for stochastic bilinear saddle-point problems with composites based on the state-of-the-art method of \cite{kovalev2021accelerated}. The latter method, for the considered class of problems, works better than Stochastic Accelerated Primal-Dual method \cite{zhang2021robust} which we use in Section~\ref{S-PP}.  
\\[2pt]
Interestingly, our \textit{Smoothing scheme} is not the only technique for reduction of algorithms for smooth convex problems to algorithms for non-smooth convex problems. \textit{Universal Nesterov's technique} from \cite{nesterov2015universal} claims that if $f$ is $M_2$-Lipschitz, then, for arbitrary $\delta > 0$,
\[ f(y) \le f(x) + \langle\nabla f(x), y-x \rangle + \frac{L}{2}\|y-x\|_2^2 + \delta,\]
with $L = M_2^2/(2\delta)$. 
Based on this observation, a generic reduction of algorithms with inexact gradient for smooth convex problems to algorithms for non-smooth problems is possible. 
\dk{Although} it may seem that this approach combined with full gradient approximation by finite-differences \cite{berahas2019theoretical} can be applied for our setting, it is not the case.
Firstly, since $L = M_2^2/(2\delta)$ and for Accelerated gradient method $\delta$ should be such that $N\delta \simeq \varepsilon$, since $\beta = 1$, we obtain that $L \sim \varepsilon^{-3/2}$, rather than $L \sim \varepsilon^{-1}$ that we have in our \textit{Smoothing scheme}. Secondly, the full gradient approximation by finite-differences provably works for smooth objective $f$ \cite{berahas2019theoretical}, which is not our case. 
\\[2pt]
Thus, the proposed \textit{Smoothing scheme} is better than alternative approaches based on other smoothing techniques such as the \textit{Universal Nesterov's technique}  \cite{nesterov2015universal}, \textit{Nesterov's smoothing via regularization of dual problem} \cite{nesterov2005smooth} and its different generalizations \cite{dev12,nguyen2017smoothing,tran2017adaptive}, which are designed and work good in the setting of first-order methods, but lead to inferior complexity in the zeroth-order setting. Unlike these approaches, our \textit{Smoothing scheme} achieves the best known complexity in terms of the successive iterations number with the best possible number of zeroth-order oracle calls. 

\section{Experiments}
\subsection{Reinforcement learning}
\begin{figure}
    \centering
    \vskip-5pt
    \includegraphics[width=0.45\textwidth]{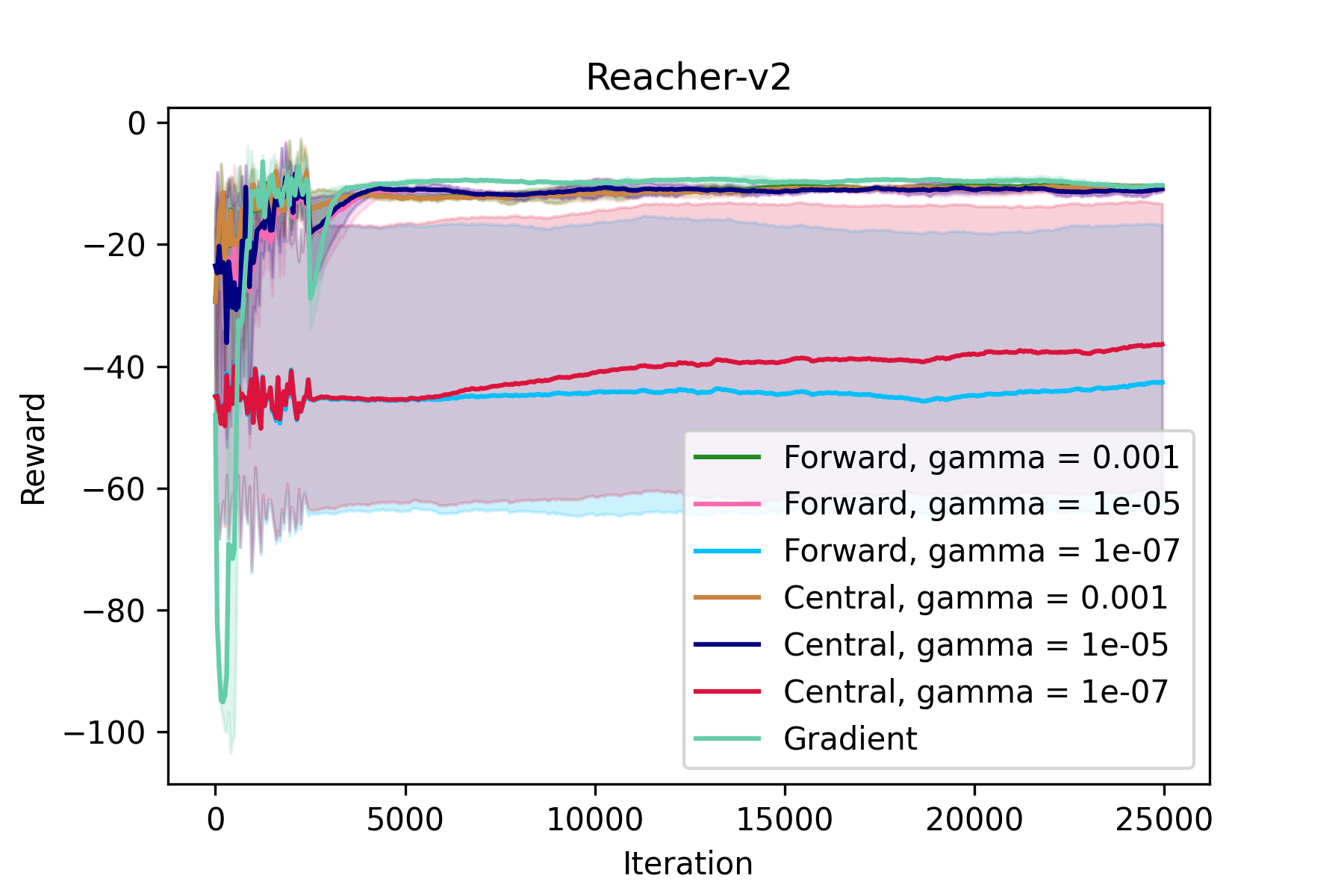}
    \vskip-10pt
    \caption{Actor's reward for ADAM with Forward and Central differences for various $\gamma$ and exact gradient ADAM. \textit{lr=} 1e-5.}
    \label{rl_main}
    \vskip-5pt
    \centering
    \vskip-5pt
    \includegraphics[width=0.45\textwidth]{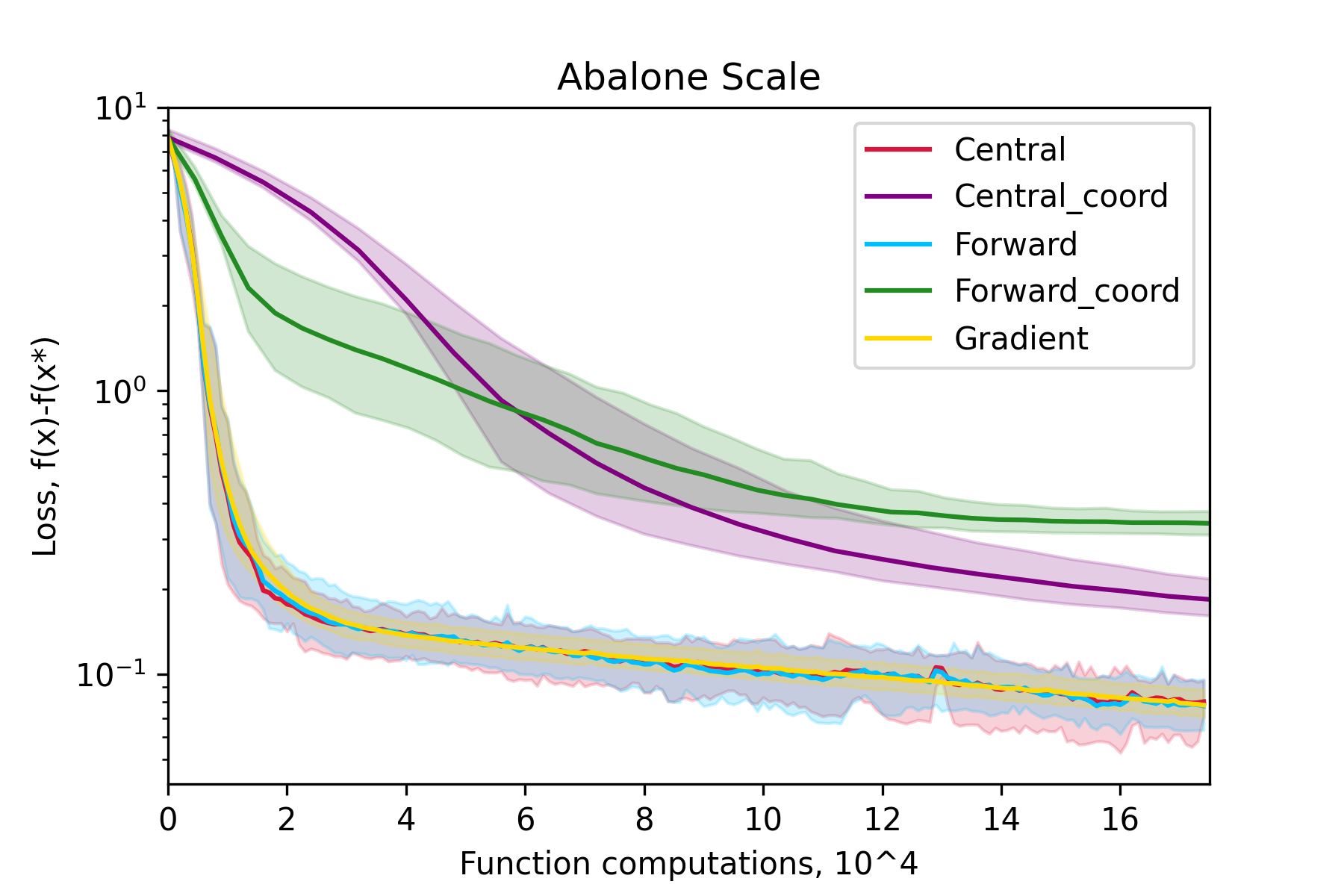}
    \vskip-10pt
    \caption{Loss for \textit{abalone scale} dataset with \textit{batch size} $=100$, learning rate is $0.1$ and $\gamma =$1e-5.}
    \vskip-10pt
    \label{l1_main}
\end{figure}

Reinforcement Learning (RL) is one of the key motivations for the proposed approach. We focus in this section on the Actor-Critic architecture and assume that the Critic is available to the Actor during training through a black-box oracle. This situation naturally motivates the application of zeroth-oder methods in which we only have access to function values. Moreover, the optimization problem is stochastic and may be convex or non-convex, smooth or non-smooth. Our RL experiments are carried out in the environment called "Reacher-v2," which is provided by the Open AI Gym toolkit. The network structure is described in Appendix \ref{appendix_RL}. We use PyTorch ADAM optimizer with three alternative inexact gradients for the Actor learning problem: exact gradient~(Gradient), central finite difference~(Central) from \eqref{sg} and
forward finite difference~(Forward) defined as
\begin{equation}\label{forward}
    \nabla f_{\gamma}(x,e) = d\frac{f(x+\gamma e) - f(x)}{\gamma}e.
\end{equation}
For a better visualization, we plot the moving average of the reward with a window size of $250$. Figure \ref{rl_main} shows the convergence of our methods.
\subsection{Robust Linear Regression}
Least absolute deviation (LAD) is a non-smooth convex problem, one of the variant of a robust linear regression \cite{yu2017robust}. It is more robust to outliers in data than the standard Linear Regression. The problem statement is 
\begin{equation}
    \min\limits_{w\in \mathbb{R}^d}  \lbrace f(w) = \tfrac{1}{n}\textstyle{\sum}_{k=1}^{m} |x_i^T w - y_i|  \rbrace,
\end{equation}
where $(x_i,y_i)$ are feature and target pairs. For the experiments we take simple dataset "abalone scale" from the LibSVM \cite{CC01a}. For this problem, we also implement ZO method with central coordinate and forward coordinate finite differences for the comparison. Their full definitions can be found in the Appendix \ref{appendix_l1}. Figure \ref{l1_main} shows the convergence of our methods.

\subsection{Support Vector Machine}
\begin{figure}
\vskip-5pt
    \centering
    {
    \includegraphics[trim={0 0 0 10pt},clip,
    width=0.45\textwidth]{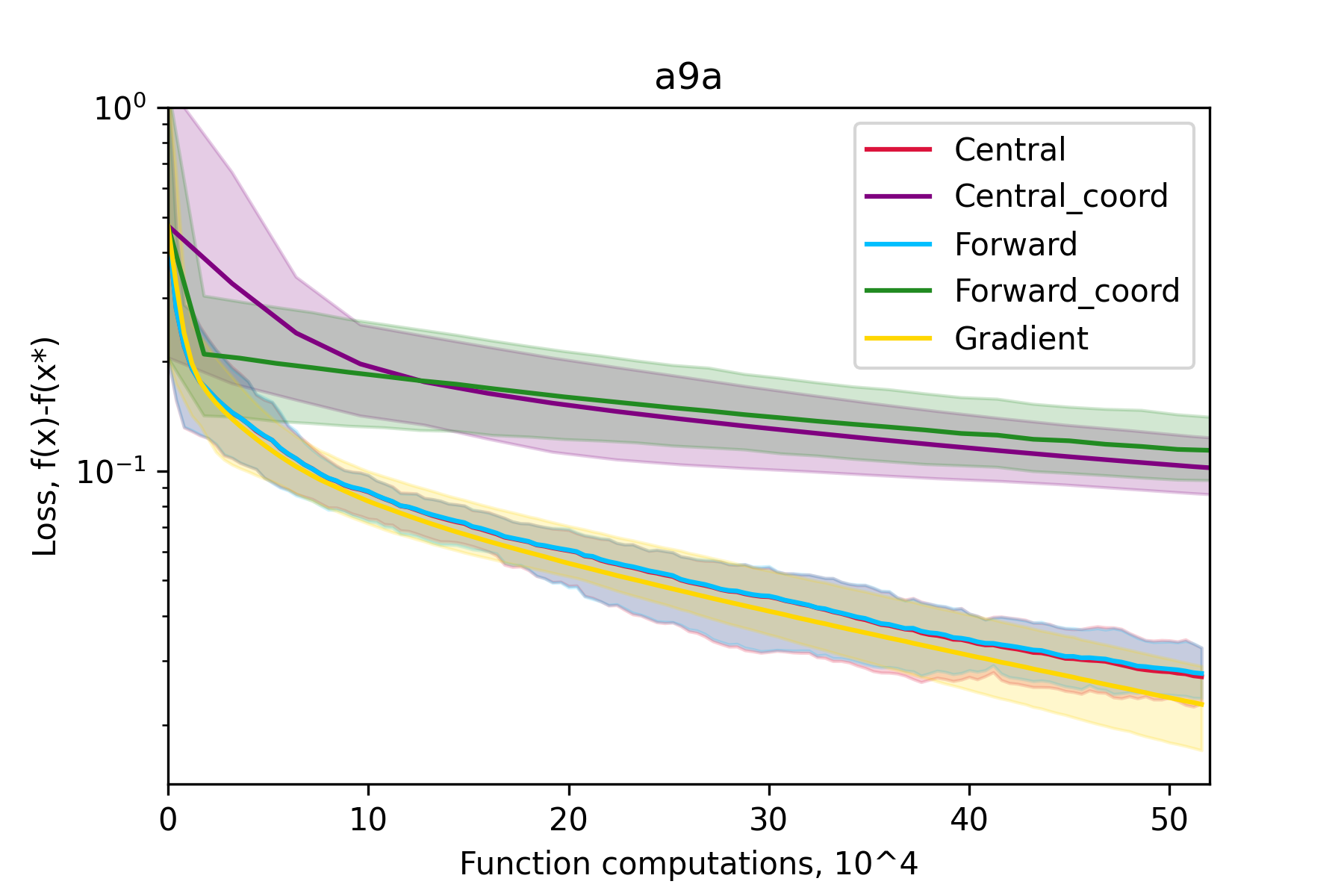}
    }
    \vskip-10pt
    \caption{Loss for \textit{a9a} dataset with $\mu=$1e-05, \textit{batch size} $=100$, $lr = 0.1$. $\gamma =$1e-5.}
    \label{svm_main}
    \vskip-10pt
\end{figure}

Support Vector Machine (SVM) is one of the classical classification algorithms that is still very popular. The problem statement is:
\begin{equation}
    \min\limits_{w\in \mathbb{R}^d}  \lbrace f(w) = \tfrac{\mu}{2}\|w\|^2 + \tfrac{1}{n}\textstyle{\sum}_{k=1}^{m} (1- y_i\cdot x_i^T w)_{+}  \rbrace,
\end{equation}
where $(x_i,y_i)$ are feature and label pairs. We use the LibSVM basic dataset "a9a" for our experiments with this problem. Figure \ref{svm_main} shows the convergence of our methods.
\subsection{Conclusion on Experiments}
We run experiments with $4$ different seeds ranging from $0$ to $3$ and present mean values as a main line and a light filling between the maximum and minimum values.
The Figures \ref{rl_main}, \ref{l1_main} and \ref{svm_main} show that with a decent choice of $\gamma$,  ADAM with zeroth-order oracle works nearly as well as  ADAM with exact gradient. Central and Forward differences outperform their coordinate versions, as it can be seen from Figures \ref{l1_main} and \ref{svm_main}. The effect of computational instability appears for tiny $\gamma$ since $x$ is too close to $x+\gamma e$ and we reach the machine precision when computing the function value. However, Figure \ref{rl_main} demonstrates that the Central approximation is more robust than Forward to such errors since the distance between points is twice larger than in the Forward case. All these plots show that our approach is applicable and can compete with gradient-based methods. For more experiments, see Appendixes \ref{appendix_RL}, \ref{appendix_l1}, \ref{appendix_svm}.
\section*{Conclusion}
This paper is devoted to the development of a universal \textit{smoothing scheme} that allows \dk{us} to construct efficient zeroth-order/gradient-free methods for non-smooth problems based on batched-gradient methods for smooth problems. This scheme preserves the  efficiency of the input methods transferring it to the output zeroth-order method. Our \textit{smoothing scheme} combines well with many other reduction techniques (batching, restarts, Catalyst acceleration for saddle-point problems, consensus-projection for decentralized distributed algorithms) that allows \dk{us} to obtain via our scheme many algorithms for a wide class of non-smooth problems. As a future work we mention adaptive 
\cite{DVINSKIKH20201715,Ene_Nguyen_Vladu_2021} and inexact model \cite{stonyakin2021} (composite, max-structure, etc.) generalizations. \ag{Another generalization consist in replacing $2$-sphere randomization on $1$-sphere randomization in smoothing scheme \cite{gasnikov2016gradient}. In some special regimes such randomization allows to improve complexity estimates on a logarithmic factor \cite{akhavan2022gradient}.}

\section*{Acknowledgements}
The work of A. Gasnikov was supported by a grant for research centers in the field of artificial intelligence, provided by the Analytical Center for the Government of the Russian Federation in accordance with the subsidy agreement (agreement identifier 000000D730321P5Q0002) and the agreement with the Ivannikov Institute for System Programming of the Russian Academy of Sciences dated November 2, 2021 No. 70-2021-00142.

\bibliography{example_paper}

\begin{thebibliography}{10}

\bibitem{agafonov2020inexact}
A.~Agafonov, D.~Kamzolov, P.~Dvurechensky, and A.~Gasnikov.
\newblock Inexact tensor methods and their application to stochastic convex
  optimization.
\newblock {\em arXiv preprint arXiv:2012.15636}, 2020.

\bibitem{akhavan2022gradient}
A.~Akhavan, E.~Chzhen, M.~Pontil, and A.~B. Tsybakov.
\newblock A gradient estimator via l1-randomization for online zero-order
  optimization with two point feedback.
\newblock {\em arXiv preprint arXiv:2205.13910}, 2022.

\bibitem{akhavan2020exploiting}
A.~Akhavan, M.~Pontil, and A.~Tsybakov.
\newblock Exploiting higher order smoothness in derivative-free optimization
  and continuous bandits.
\newblock {\em Advances in Neural Information Processing Systems},
  33:9017--9027, 2020.

\bibitem{bayandina2018gradient}
A.~S. Bayandina, A.~V. Gasnikov, and A.~A. Lagunovskaya.
\newblock Gradient-free two-point methods for solving stochastic nonsmooth
  convex optimization problems with small non-random noises.
\newblock {\em Automation and Remote Control}, 79(8):1399--1408, 2018.

\bibitem{berahas2019theoretical}
A.~S. Berahas, L.~Cao, K.~Choromanski, and K.~Scheinberg.
\newblock A theoretical and empirical comparison of gradient approximations in
  derivative-free optimization.
\newblock {\em arXiv:1905.01332}, 2019.

\bibitem{beznosikov2020derivative}
A.~Beznosikov, E.~Gorbunov, and A.~Gasnikov.
\newblock Derivative-free method for composite optimization with applications
  to decentralized distributed optimization.
\newblock {\em IFAC-PapersOnLine}, 53(2):4038--4043, 2020.

\bibitem{beznosikov2021one}
A.~Beznosikov, V.~Novitskii, and A.~Gasnikov.
\newblock One-point gradient-free methods for smooth and non-smooth
  saddle-point problems.
\newblock In {\em International Conference on Mathematical Optimization Theory
  and Operations Research}, pages 144--158. Springer, 2021.

\bibitem{beznosikov2020gradient}
A.~Beznosikov, A.~Sadiev, and A.~Gasnikov.
\newblock Gradient-free methods with inexact oracle for convex-concave
  stochastic saddle-point problem.
\newblock In {\em International Conference on Mathematical Optimization Theory
  and Operations Research}, pages 105--119. Springer, 2020.

\bibitem{aleks2020distributed}
A.~Beznosikov, V.~Samokhin, and A.~Gasnikov.
\newblock Distributed saddle-point problems: Lower bounds, optimal and robust
  algorithms, 2020.

\bibitem{bubeck2019complexity}
S.~Bubeck, Q.~Jiang, Y.~T. Lee, Y.~Li, A.~Sidford, et~al.
\newblock Complexity of highly parallel non-smooth convex optimization.
\newblock {\em Advances in neural information processing systems}, 2019.

\bibitem{CC01a}
C.-C. Chang and C.-J. Lin.
\newblock {LIBSVM}: A library for support vector machines.
\newblock {\em ACM Transactions on Intelligent Systems and Technology},
  2:27:1--27:27, 2011.
\newblock Software available at \url{http://www.csie.ntu.edu.tw/~cjlin/libsvm}.

\bibitem{chen2017}
P.-Y. Chen, H.~Zhang, Y.~Sharma, J.~Yi, and C.-J. Hsieh.
\newblock Zoo.
\newblock {\em Proceedings of the 10th ACM Workshop on Artificial Intelligence
  and Security}, Nov 2017.

\bibitem{conn2009introduction}
A.~R. Conn, K.~Scheinberg, and L.~N. Vicente.
\newblock {\em Introduction to Derivative-Free Optimization}.
\newblock Society for Industrial and Applied Mathematics, 2009.

\bibitem{cotter2011better}
A.~Cotter, O.~Shamir, N.~Srebro, and K.~Sridharan.
\newblock Better mini-batch algorithms via accelerated gradient methods.
\newblock {\em Advances in Neural Information Processing Systems},
  24:1647--1655, 2011.

\bibitem{devolder2013exactness}
O.~Devolder.
\newblock {\em Exactness, inexactness and stochasticity in first-order methods
  for large-scale convex optimization}.
\newblock PhD thesis, PhD thesis, 2013.

\bibitem{dev12}
O.~Devolder, F.~Glineur, and Y.~Nesterov.
\newblock Double smoothing technique for large-scale linearly constrained
  convex optimization.
\newblock {\em SIAM Journal on Optimization}, 22(2):702--727, 2012.

\bibitem{diakonikolas2020lower}
J.~Diakonikolas and C.~Guzm{\'a}n.
\newblock Lower bounds for parallel and randomized convex optimization.
\newblock {\em J. Mach. Learn. Res.}, 21:5--1, 2020.

\bibitem{duchi2012randomized}
J.~C. Duchi, P.~L. Bartlett, and M.~J. Wainwright.
\newblock Randomized smoothing for stochastic optimization.
\newblock {\em SIAM Journal on Optimization}, 22(2):674--701, 2012.

\bibitem{duchi2015optimal}
J.~C. Duchi, M.~I. Jordan, M.~J. Wainwright, and A.~Wibisono.
\newblock Optimal rates for zero-order convex optimization: The power of two
  function evaluations.
\newblock {\em IEEE Transactions on Information Theory}, 61(5):2788--2806,
  2015.

\bibitem{Dvinskikh2021}
D.~Dvinskikh and A.~Gasnikov.
\newblock Decentralized and parallel primal and dual accelerated methods for
  stochastic convex programming problems.
\newblock {\em Journal of Inverse and Ill-posed Problems}, 29(3):385--405,
  2021.

\bibitem{DVINSKIKH20201715}
D.~Dvinskikh, A.~Ogaltsov, A.~Gasnikov, P.~Dvurechensky, and V.~Spokoiny.
\newblock On the line-search gradient methods for stochastic optimization.
\newblock {\em IFAC-PapersOnLine}, 53(2):1715--1720, 2020.
\newblock 21st IFAC World Congress.

\bibitem{dvinskikh2022gradient}
D.~Dvinskikh, V.~Tominin, Y.~Tominin, and A.~Gasnikov.
\newblock Gradient-free optimization for non-smooth minimax problems with
  maximum value of adversarial noise.
\newblock {\em arXiv preprint arXiv:2202.06114}, 2022.

\bibitem{dvinskikh2020accelerated}
D.~Dvinskikh, A.~Tyurin, A.~Gasnikov, and S.~Omelchenko.
\newblock Accelerated and nonaccelerated stochastic gradient descent with model
  conception.
\newblock {\em Math. Notes}, 108(4):511--522, 2020.

\bibitem{dvurechensky2016stochastic}
P.~Dvurechensky and A.~Gasnikov.
\newblock Stochastic intermediate gradient method for convex problems with
  stochastic inexact oracle.
\newblock {\em Journal of Optimization Theory and Applications},
  171(1):121--145, 2016.

\bibitem{dvurechensky2021accelerated}
P.~Dvurechensky, E.~Gorbunov, and A.~Gasnikov.
\newblock An accelerated directional derivative method for smooth stochastic
  convex optimization.
\newblock {\em European Journal of Operational Research}, 290(2):601--621,
  2021.

\bibitem{Ene_Nguyen_Vladu_2021}
A.~Ene, H.~L. Nguyen, and A.~Vladu.
\newblock Adaptive gradient methods for constrained convex optimization and
  variational inequalities.
\newblock {\em Proceedings of the AAAI Conference on Artificial Intelligence},
  35(8):7314--7321, May 2021.

\bibitem{ermoliev1976stochastic}
Y.~Ermoliev.
\newblock Stochastic programming methods, 1976.

\bibitem{feldman2016generalization}
V.~Feldman.
\newblock Generalization of erm in stochastic convex optimization: The
  dimension strikes back.
\newblock {\em Advances in Neural Information Processing Systems},
  29:3576--3584, 2016.

\bibitem{flaxman2005online}
A.~D. Flaxman, A.~T. Kalai, and H.~B. McMahan.
\newblock Online convex optimization in the bandit setting: gradient descent
  without a gradient.
\newblock In {\em Proceedings of the sixteenth annual ACM-SIAM symposium on
  Discrete algorithms}, pages 385--394, 2005.

\bibitem{gasnikov2017universal}
A.~Gasnikov.
\newblock Universal gradient descent.
\newblock {\em MCCME, arXiv:1711.00394}, 2021.

\bibitem{gasnikov2019near}
A.~Gasnikov, P.~Dvurechensky, E.~Gorbunov, E.~Vorontsova, D.~Selikhanovych,
  C.~A. Uribe, B.~Jiang, H.~Wang, S.~Zhang, S.~Bubeck, et~al.
\newblock Near optimal methods for minimizing convex functions with lipschitz $
  p $-th derivatives.
\newblock In {\em Conference on Learning Theory}, pages 1392--1393. PMLR, 2019.

\bibitem{gasnikov2017stochastic}
A.~V. Gasnikov, E.~A. Krymova, A.~A. Lagunovskaya, I.~N. Usmanova, and F.~A.
  Fedorenko.
\newblock Stochastic online optimization. single-point and multi-point
  non-linear multi-armed bandits. convex and strongly-convex case.
\newblock {\em Automation and remote control}, 78(2):224--234, 2017.

\bibitem{gasnikov2016gradient}
A.~V. Gasnikov, A.~A. Lagunovskaya, I.~N. Usmanova, and F.~A. Fedorenko.
\newblock Gradient-free proximal methods with inexact oracle for convex
  stochastic nonsmooth optimization problems on the simplex.
\newblock {\em Automation and Remote Control}, 77(11):2018--2034, 2016.

\bibitem{gladin2021solving}
E.~Gladin, A.~Sadiev, A.~Gasnikov, P.~Dvurechensky, A.~Beznosikov, and
  M.~Alkousa.
\newblock Solving smooth min-min and min-max problems by mixed oracle
  algorithms.
\newblock In {\em International Conference on Mathematical Optimization Theory
  and Operations Research}. Springer, 2021.

\bibitem{gorbunov2021stochastic}
E.~Gorbunov, H.~Berard, G.~Gidel, and N.~Loizou.
\newblock Stochastic extragradient: General analysis and improved rates, 2021.

\bibitem{gorbunov2020stochastic}
E.~Gorbunov, M.~Danilova, and A.~Gasnikov.
\newblock Stochastic optimization with heavy-tailed noise via accelerated
  gradient clipping.
\newblock {\em Advances in Neural Information Processing Systems},
  33:15042--15053, 2020.

\bibitem{gorbunov2019optimal}
E.~Gorbunov, D.~Dvinskikh, and A.~Gasnikov.
\newblock Optimal decentralized distributed algorithms for stochastic convex
  optimization.
\newblock {\em arXiv preprint arXiv:1911.07363}, 2019.

\bibitem{gorbunov2018accelerated}
E.~Gorbunov, P.~Dvurechensky, and A.~Gasnikov.
\newblock An accelerated method for derivative-free smooth stochastic convex
  optimization.
\newblock {\em SIAM J. Optim. arXiv:1802.09022}, 2022.

\bibitem{gorbunov2020recent}
E.~Gorbunov, A.~Rogozin, A.~Beznosikov, D.~Dvinskikh, and A.~Gasnikov.
\newblock Recent theoretical advances in decentralized distributed convex
  optimization.
\newblock {\em arXiv preprint arXiv:2011.13259}, 2020.

\bibitem{gorbunov2019upper}
E.~Gorbunov, E.~A. Vorontsova, and A.~V. Gasnikov.
\newblock On the upper bound for the expectation of the norm of a vector
  uniformly distributed on the sphere and the phenomenon of concentration of
  uniform measure on the sphere.
\newblock {\em Mathematical Notes}, 106, 2019.

\bibitem{juditsky2011solving}
A.~Juditsky, A.~Nemirovski, and C.~Tauvel.
\newblock Solving variational inequalities with stochastic mirror-prox
  algorithm.
\newblock {\em Stochastic Systems}, 1(1):17--58, 2011.

\bibitem{juditsky2008large}
A.~Juditsky and A.~S. Nemirovski.
\newblock Large deviations of vector-valued martingales in 2-smooth normed
  spaces.
\newblock {\em arXiv preprint arXiv:0809.0813}, 2008.

\bibitem{juditsky2014deterministic}
A.~Juditsky and Y.~Nesterov.
\newblock Deterministic and stochastic primal-dual subgradient algorithms for
  uniformly convex minimization.
\newblock {\em Stochastic Systems}, 4(1):44--80, 2014.

\bibitem{kotsalis2020simple}
G.~Kotsalis, G.~Lan, and T.~Li.
\newblock Simple and optimal methods for stochastic variational inequalities,
  i: operator extrapolation.
\newblock {\em arXiv preprint arXiv:2011.02987}, 2020.

\bibitem{kovalev2021accelerated}
D.~Kovalev, A.~Gasnikov, and P.~Richtárik.
\newblock Accelerated primal-dual gradient method for smooth and convex-concave
  saddle-point problems with bilinear coupling, 2021.

\bibitem{kulunchakov2020optimisation}
A.~Kulunchakov.
\newblock {\em Optimisation stochastique pour l'apprentissage machine {\`a}
  grande {\'e}chelle: r{\'e}duction de la variance et acc{\'e}l{\'e}ration}.
\newblock PhD thesis, Universit{\'e} Grenoble Alpes, 2020.

\bibitem{kulunchakov2019estimate}
A.~Kulunchakov and J.~Mairal.
\newblock Estimate sequences for variance-reduced stochastic composite
  optimization.
\newblock In {\em International Conference on Machine Learning}, pages
  3541--3550. PMLR, 2019.

\bibitem{lan2012optimal}
G.~Lan.
\newblock An optimal method for stochastic composite optimization.
\newblock {\em Mathematical Programming}, 133(1):365--397, 2012.

\bibitem{lan2020first}
G.~Lan.
\newblock {\em First-order and Stochastic Optimization Methods for Machine
  Learning}.
\newblock Springer, 2020.

\bibitem{lan2018random}
G.~Lan and Y.~Zhou.
\newblock Random gradient extrapolation for distributed and stochastic
  optimization.
\newblock {\em SIAM Journal on Optimization}, 28(4):2753--2782, 2018.

\bibitem{lillicrap2019continuous}
T.~P. Lillicrap, J.~J. Hunt, A.~Pritzel, N.~Heess, T.~Erez, Y.~Tassa,
  D.~Silver, and D.~Wierstra.
\newblock Continuous control with deep reinforcement learning, 2019.

\bibitem{pmlr-v125-lin20a}
T.~Lin, C.~Jin, and M.~I. Jordan.
\newblock Near-optimal algorithms for minimax optimization.
\newblock In J.~Abernethy and S.~Agarwal, editors, {\em Proceedings of Thirty
  Third Conference on Learning Theory}, volume 125 of {\em Proceedings of
  Machine Learning Research}, pages 2738--2779. PMLR, 09--12 Jul 2020.

\bibitem{nemirovsky1983problem}
A.~Nemirovsky and D.~Yudin.
\newblock Problem complexity and method efficiency in optimization.-j. wiley \&
  sons, new york.
\newblock 1983.

\bibitem{nesterov2005smooth}
Y.~Nesterov.
\newblock Smooth minimization of non-smooth functions.
\newblock {\em Mathematical Programming}, 103(1):127--152, 2005.

\bibitem{nesterov2015universal}
Y.~Nesterov.
\newblock Universal gradient methods for convex optimization problems.
\newblock {\em Mathematical Programming}, 152(1):381--404, 2015.

\bibitem{nesterov2021implementable}
Y.~Nesterov.
\newblock Implementable tensor methods in unconstrained convex optimization.
\newblock {\em Mathematical Programming}, 186(1):157--183, 2021.

\bibitem{nesterov2017random}
Y.~Nesterov and V.~Spokoiny.
\newblock Random gradient-free minimization of convex functions.
\newblock {\em Foundations of Computational Mathematics}, 17(2):527--566, 2017.

\bibitem{nguyen2017smoothing}
Q.~V. Nguyen, O.~Fercoq, and V.~Cevher.
\newblock Smoothing technique for nonsmooth composite minimization with linear
  operator.
\newblock {\em arXiv:1706.05837}, 2017.

\bibitem{peters2008natural}
J.~Peters and S.~Schaal.
\newblock Natural actor-critic.
\newblock {\em Neurocomputing}, 71(7-9):1180--1190, 2008.

\bibitem{polyak1987introduction}
B.~T. Polyak.
\newblock Introduction to optimization. optimization software.
\newblock {\em Inc., Publications Division, New York}, 1, 1987.

\bibitem{risteski2016algorithms}
A.~Risteski and Y.~Li.
\newblock Algorithms and matching lower bounds for approximately-convex
  optimization.
\newblock {\em Advances in Neural Information Processing Systems},
  29:4745--4753, 2016.

\bibitem{rogozin2021accelerated}
A.~Rogozin, M.~Bochko, P.~Dvurechensky, A.~Gasnikov, and V.~Lukoshkin.
\newblock An accelerated method for decentralized distributed stochastic
  optimization over time-varying graphs.
\newblock {\em 2021 IEEE Conference on Decision and Control (CDC).
  arXiv:2103.15598}, 2021.

\bibitem{sadiev2021zeroth}
A.~Sadiev, A.~Beznosikov, P.~Dvurechensky, and A.~Gasnikov.
\newblock Zeroth-order algorithms for smooth saddle-point problems.
\newblock In {\em International Conference on Mathematical Optimization Theory
  and Operations Research}, pages 71--85. Springer, 2021.

\bibitem{scaman2019optimal}
K.~Scaman, F.~Bach, S.~Bubeck, Y.~Lee, and L.~Massouli{\'e}.
\newblock Optimal convergence rates for convex distributed optimization in
  networks.
\newblock {\em Journal of Machine Learning Research}, 20:1--31, 2019.

\bibitem{shalev2014understanding}
S.~Shalev-Shwartz and S.~Ben-David.
\newblock {\em Understanding machine learning: From theory to algorithms}.
\newblock Cambridge university press, 2014.

\bibitem{shalev2009stochastic}
S.~Shalev-Shwartz, O.~Shamir, N.~Srebro, and K.~Sridharan.
\newblock Stochastic convex optimization.
\newblock In {\em COLT}, volume~2, page~5, 2009.

\bibitem{shamir2017optimal}
O.~Shamir.
\newblock An optimal algorithm for bandit and zero-order convex optimization
  with two-point feedback.
\newblock {\em The Journal of Machine Learning Research}, 18(1):1703--1713,
  2017.

\bibitem{shapiro2021lectures}
A.~Shapiro, D.~Dentcheva, and A.~Ruszczynski.
\newblock {\em Lectures on stochastic programming: modeling and theory}.
\newblock SIAM, 2021.

\bibitem{shapiro2005complexity}
A.~Shapiro and A.~Nemirovski.
\newblock On complexity of stochastic programming problems.
\newblock In {\em Continuous optimization}, pages 111--146. Springer, 2005.

\bibitem{spall2005introduction}
J.~C. Spall.
\newblock {\em Introduction to stochastic search and optimization: estimation,
  simulation, and control}, volume~65.
\newblock John Wiley \& Sons, 2005.

\bibitem{stepanov2021onepoint}
I.~Stepanov, A.~Voronov, A.~Beznosikov, and A.~Gasnikov.
\newblock One-point gradient-free methods for composite optimization with
  applications to distributed optimization, 2021.

\bibitem{stonyakin2021}
F.~Stonyakin, A.~Tyurin, A.~Gasnikov, P.~Dvurechensky, A.~Agafonov,
  D.~Dvinskikh, M.~Alkousa, D.~Pasechnyuk, S.~Artamonov, and V.~Piskunova.
\newblock Inexact model: a framework for optimization and variational
  inequalities.
\newblock {\em Optimization Methods and Software}, 0(0):1--47, 2021.

\bibitem{cite-key}
A.~B. Taylor, J.~M. Hendrickx, and F.~Glineur.
\newblock Smooth strongly convex interpolation and exact worst-case performance
  of first-order methods.
\newblock {\em Mathematical Programming}, 161(1):307--345, 2017.

\bibitem{tian2021acceleration}
Y.~Tian, G.~Scutari, T.~Cao, and A.~Gasnikov.
\newblock Acceleration in distributed optimization under similarity, 2021.

\bibitem{vladislav2021accelerated}
V.~Tominin, Y.~Tominin, E.~Borodich, D.~Kovalev, A.~Gasnikov, and
  P.~Dvurechensky.
\newblock On accelerated methods for saddle-point problems with composite
  structure.
\newblock {\em arXiv preprint arXiv:2103.09344}, 2021.

\bibitem{tran2017adaptive}
Q.~Tran-Dinh.
\newblock Adaptive smoothing algorithms for nonsmooth composite convex
  minimization.
\newblock {\em Computational Optimization and Applications}, 66(3):425--451,
  2017.

\bibitem{woodworth2021even}
B.~Woodworth and N.~Srebro.
\newblock An even more optimal stochastic optimization algorithm: Minibatching
  and interpolation learning.
\newblock {\em arXiv preprint arXiv:2106.02720}, 2021.

\bibitem{YOUSEFIAN201256}
F.~Yousefian, A.~Nedić, and U.~V. Shanbhag.
\newblock On stochastic gradient and subgradient methods with adaptive
  steplength sequences.
\newblock {\em Automatica}, 48(1):56--67, 2012.

\bibitem{yu2017robust}
C.~Yu and W.~Yao.
\newblock Robust linear regression: A review and comparison.
\newblock {\em Communications in Statistics-Simulation and Computation},
  46(8):6261--6282, 2017.

\bibitem{zhang2021lower}
J.~Zhang, M.~Hong, and S.~Zhang.
\newblock On lower iteration complexity bounds for the convex concave saddle
  point problems.
\newblock {\em Mathematical Programming}, pages 1--35, 2021.

\bibitem{zhang2021robust}
X.~Zhang, N.~S. Aybat, and M.~Gurbuzbalaban.
\newblock Robust accelerated primal-dual methods for computing saddle points.
\newblock {\em arXiv preprint arXiv:2111.12743}, 2021.

\end{thebibliography}

\clearpage
\appendix
\onecolumn
\section{Appendix}

\subsection{Additional Experiments: Adversarial Attack}
Adversarial attack is aimed to generate input images with unobtrusive interference introduced to them to deceive a well trained classifier. These adversarial examples are created to recognize the robustness of models. In our experiments, we used ZO methods to generate adversarial examples targeted at a black-box model trained to solve a task of classification of MNIST dataset of handwritten digits. The output of this model is $F(\cdot) = [F_1(\cdot), \dots , F_K(\cdot)]$, where $F_k(\cdot)$ is a prediction score of the $k^{th}$ class. A correctly classified sample image $\mathbf{a_i}$ is taken from the dataset and its adversarial example is produced as follows: 

\begin{equation}
    \mathbf{a_i^{adv}} = 0.5\tanh (\tanh^{-1}(2\mathbf{a_i}) + \mathbf{x} )
\end{equation}

Next, we apply the same individual attacking loss utilized in \cite{chen2017}:

\begin{equation} \label{eq:aa_loss}
    f_i(\mathbf{x}) = \max \{ \log F_{y_i}(\mathbf{a_i^{adv}}) - \max_{t \neq y_i} \log F_t(\mathbf{a_i^{adv}}), 0 \} + \lambda \|\mathbf{a_i^{adv}} - \mathbf{a_i}  \|^2
\end{equation}

In our experiments we used sample images for digit class "1" and "4" and set their regularization parameters as $\lambda=1$ and $\lambda=0.1$ respectively. 

\begin{figure}
\centering
\begin{subfigure}{} 
	   \includegraphics[width=0.45\textwidth]{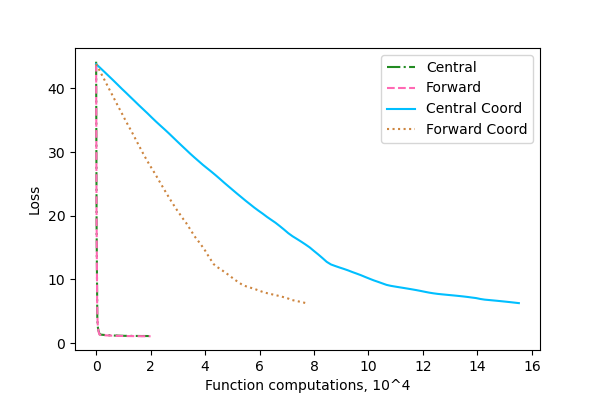}
\end{subfigure}
\hfill 
\begin{subfigure}{}
	   \includegraphics[width=0.45\textwidth]{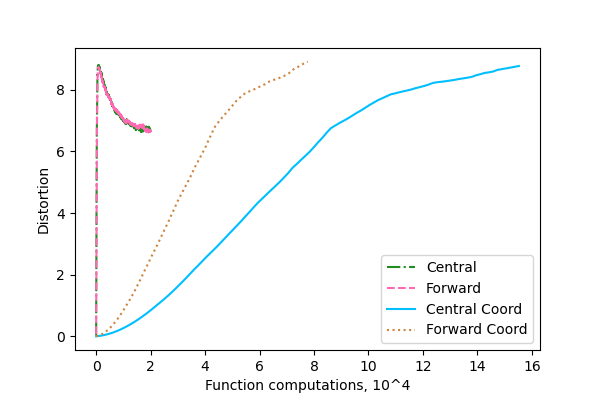} 
\end{subfigure}
\caption{Comparison of different zeroth-order algorithms for generating $n=50$ adversarial examples for digit class "4" with $\lambda=0.1$. Left: Loss (\ref{eq:aa_loss}) with Central, Forward, Central Coord and Forward Coord where $lr=0.01$ and $\gamma=0.01$. Right: Distortion($\frac{1}{n}\sum_{i=1}^n\|\mathbf{a_i^{adv}} - \mathbf{a_i}  \|^2$) average in $n=50$ generated adversarial examples.}
\label{ap:fig_adv_attack}
\end{figure}

Figure \ref{ap:fig_adv_attack} shows that methods with coordinate steps require more function computations to converge which is justified by large $d$. This can also be observed in Table \ref{table:digit1} where randomly generated adversarial examples were classified similarly by target model, although it is clear that each method produces different result.

All the experiments were conducted in Python 3 and PyTorch 1.10.1 on an Ubuntu 20.04.3 LTS machine with Intel(R) Xeon(R) Silver 4215 CPU @ 2.50GHz and 125 GB RAM.

\begin{table*}[ht]
	\small
	\caption{Generated adversarial examples for digit ``1''  class from a random batch of $n=10$ images, where image distortion is defined as $\frac{1}{n}\sum_{i=1}^n \|\mathbf{a}^{adv}_i-\mathbf{a}_i\|^2$.}
	\vspace{0.2cm}
	\begin{adjustbox}{max width=\textwidth }
		\begin{tabular}
			{cccccccccccc}
			\hline
			Image ID & $7$& $10$ & $13$ & $18$ & $19$ & $20$ & $21$ & $22$ & $28$ & $29$ & Image  distortion \\
			\hline &&&&&&&&&& \vspace{-0.2cm} \\
			$\bf Original$ &
			\parbox[c]{3em}{\includegraphics[width=0.5in]{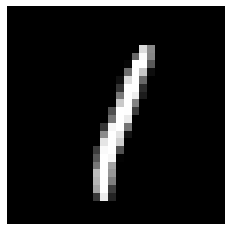}} &
			\parbox[c]{3em}{\includegraphics[width=0.5in]{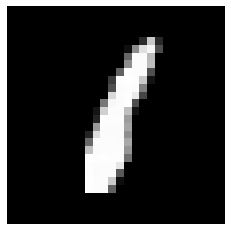}} &
			\parbox[c]{3em}{\includegraphics[width=0.5in]{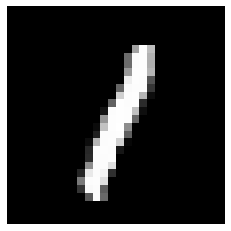}} &
			\parbox[c]{3em}{\includegraphics[width=0.5in]{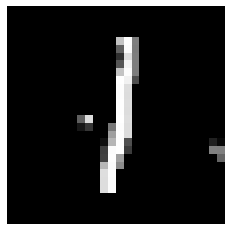}} &
			\parbox[c]{3em}{\includegraphics[width=0.5in]{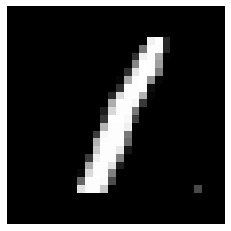}} &
			\parbox[c]{3em}{\includegraphics[width=0.5in]{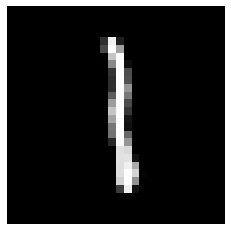}} &
			\parbox[c]{3em}{\includegraphics[width=0.5in]{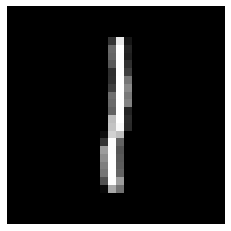}} &
			\parbox[c]{3em}{\includegraphics[width=0.5in]{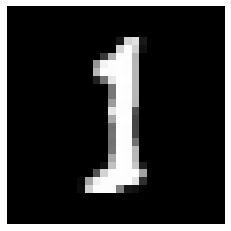}} &
			\parbox[c]{3em}{\includegraphics[width=0.5in]{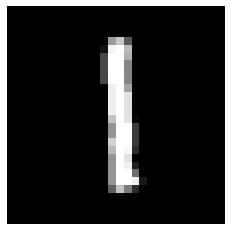}} &
			\parbox[c]{3em}{\includegraphics[width=0.5in]{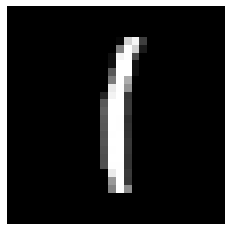}} &
			\vspace{0.1cm}\\
			
			\hline &&&&&&&&&& \vspace{-0.2cm} \\
			Central &
			\parbox[c]{3em}{\includegraphics[width=0.5in]{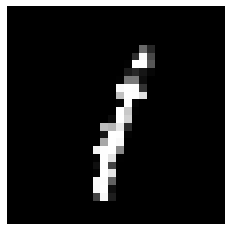}} &
			\parbox[c]{3em}{\includegraphics[width=0.5in]{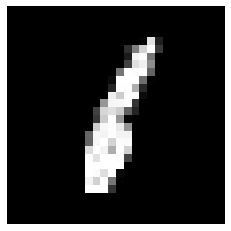}} &
			\parbox[c]{3em}{\includegraphics[width=0.5in]{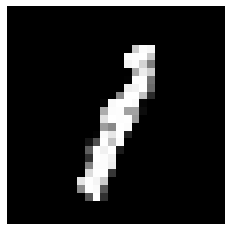}} &
			\parbox[c]{3em}{\includegraphics[width=0.5in]{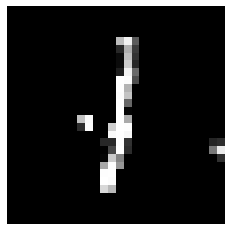}} &
			\parbox[c]{3em}{\includegraphics[width=0.5in]{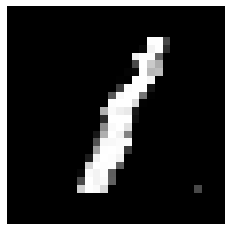}} &
			\parbox[c]{3em}{\includegraphics[width=0.5in]{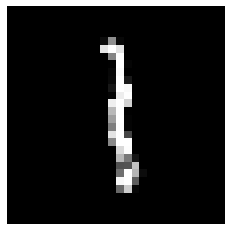}} &
			\parbox[c]{3em}{\includegraphics[width=0.5in]{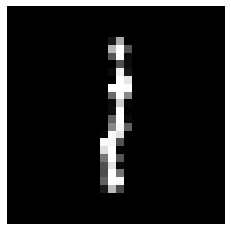}} &
			\parbox[c]{3em}{\includegraphics[width=0.5in]{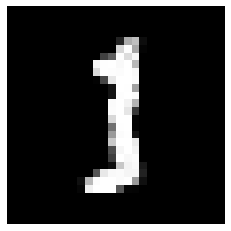}} &
			\parbox[c]{3em}{\includegraphics[width=0.5in]{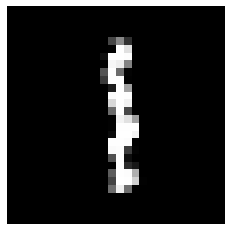}} &
			\parbox[c]{3em}{\includegraphics[width=0.5in]{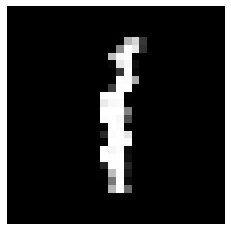}} &$5.0107$\vspace{0.1cm}\\
			\hline
			Classified as & $\bf 4$ & $\bf 6$ & $\bf 7$ & $\bf 4$ & $\bf 8$ & $\bf 8$ & $\bf 8$ & $\bf 3$ & $\bf 8$ & $\bf 8$& \\
			\hline &&&&&&&&&& \vspace{-0.2cm} \\
			
			Forward &
			\parbox[c]{3em}{\includegraphics[width=0.5in]{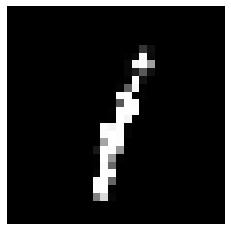}} &
			\parbox[c]{3em}{\includegraphics[width=0.5in]{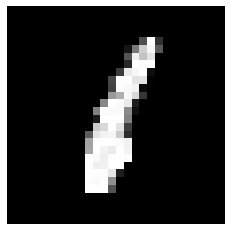}} &
			\parbox[c]{3em}{\includegraphics[width=0.5in]{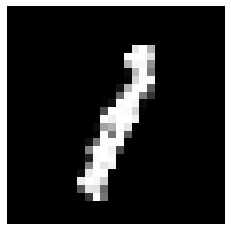}} &
			\parbox[c]{3em}{\includegraphics[width=0.5in]{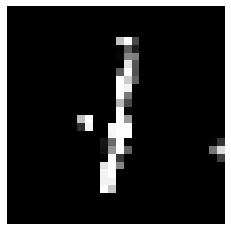}} &
			\parbox[c]{3em}{\includegraphics[width=0.5in]{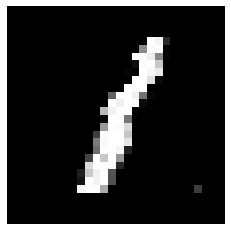}} &
			\parbox[c]{3em}{\includegraphics[width=0.5in]{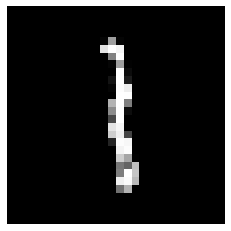}} &
			\parbox[c]{3em}{\includegraphics[width=0.5in]{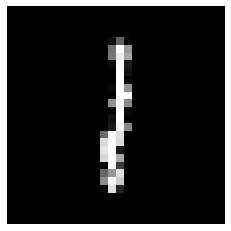}} &
			\parbox[c]{3em}{\includegraphics[width=0.5in]{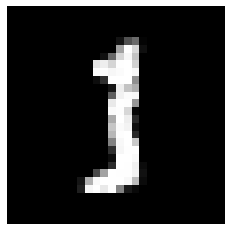}} &
			\parbox[c]{3em}{\includegraphics[width=0.5in]{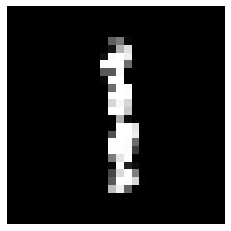}} &
			\parbox[c]{3em}{\includegraphics[width=0.5in]{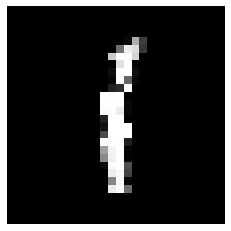}} &$5.1004$\vspace{0.1cm}\\
			\hline
			Classified as & $\bf 4$ & $\bf 6$ & $\bf 7$ & $\bf 4$ & $\bf 8$ & $\bf 8$ & $\bf 8$ & $\bf 3$ & $\bf 8$ & $\bf 8$& \\
			\hline &&&&&&&&&& \vspace{-0.2cm} \\
			
			Central Coord &
			\parbox[c]{3em}{\includegraphics[width=0.5in]{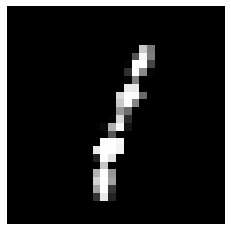}} &
			\parbox[c]{3em}{\includegraphics[width=0.5in]{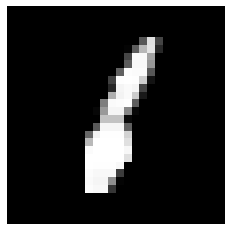}} &
			\parbox[c]{3em}{\includegraphics[width=0.5in]{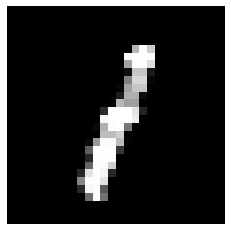}} &
			\parbox[c]{3em}{\includegraphics[width=0.5in]{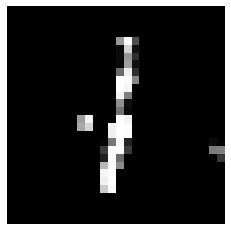}} &
			\parbox[c]{3em}{\includegraphics[width=0.5in]{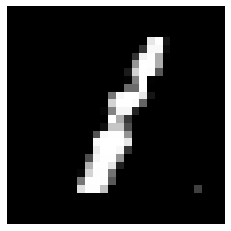}} &
			\parbox[c]{3em}{\includegraphics[width=0.5in]{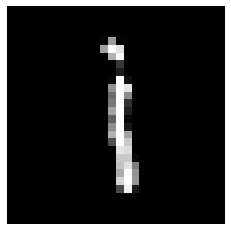}} &
			\parbox[c]{3em}{\includegraphics[width=0.5in]{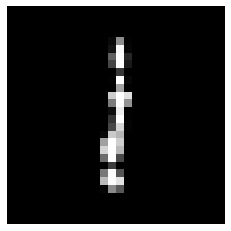}} &
			\parbox[c]{3em}{\includegraphics[width=0.5in]{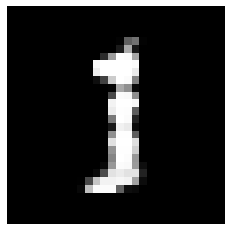}} &
			\parbox[c]{3em}{\includegraphics[width=0.5in]{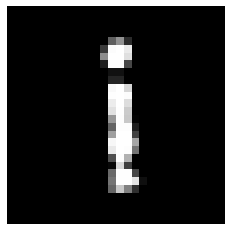}} &
			\parbox[c]{3em}{\includegraphics[width=0.5in]{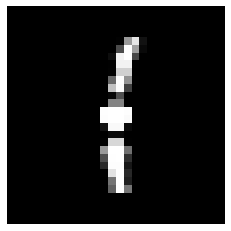}} &$4.6761$\vspace{0.1cm}\\
			\hline
			Classified as & $\bf 4$ & $\bf 6$ & $\bf 7$ & $\bf 4$ & $\bf 8$ & $\bf 8$ & $\bf 8$ & $\bf 3$ & $\bf 8$ & $\bf 4$& \\
			\hline &&&&&&&&&& \vspace{-0.2cm} \\
			
			Forward Coord &
			\parbox[c]{3em}{\includegraphics[width=0.5in]{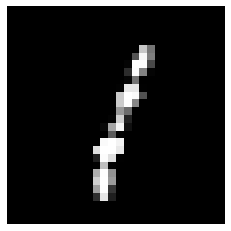}} &
			\parbox[c]{3em}{\includegraphics[width=0.5in]{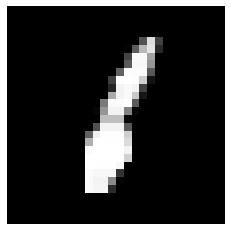}} &
			\parbox[c]{3em}{\includegraphics[width=0.5in]{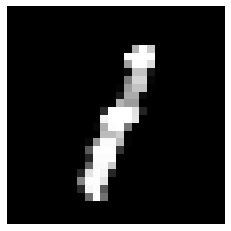}} &
			\parbox[c]{3em}{\includegraphics[width=0.5in]{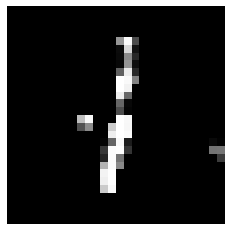}} &
			\parbox[c]{3em}{\includegraphics[width=0.5in]{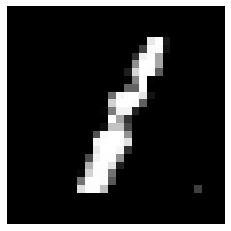}} &
			\parbox[c]{3em}{\includegraphics[width=0.5in]{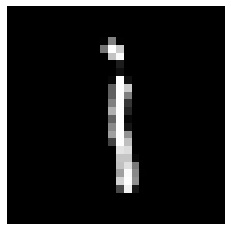}} &
			\parbox[c]{3em}{\includegraphics[width=0.5in]{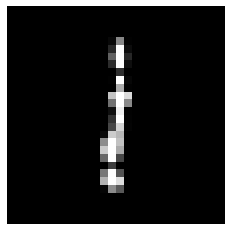}} &
			\parbox[c]{3em}{\includegraphics[width=0.5in]{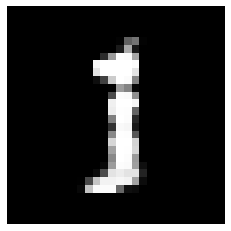}} &
			\parbox[c]{3em}{\includegraphics[width=0.5in]{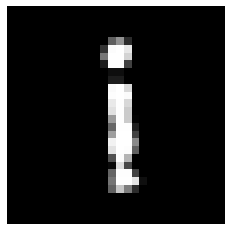}} &
			\parbox[c]{3em}{\includegraphics[width=0.5in]{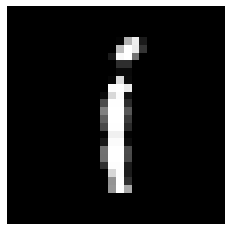}} &$4.4163$\vspace{0.1cm}\\
			\hline
			Classified as & $\bf 4$ & $\bf 6$ & $\bf 7$ & $\bf 4$ & $\bf 8$ & $\bf 8$ & $\bf 8$ & $\bf 3$ & $\bf 8$ & $\bf 8$& \\
			\hline

		\end{tabular}
	\end{adjustbox}
	\label{table:digit1}
\end{table*}

\subsection{Additional Experiments: RL experiments}
\label{appendix_RL}
Our RL experiments are carried out in an environment called "Reacher-v2," which is provided by the Open AI Gym toolkit. This environment simulates an agent ($2$ DOF robotic arm) tasked with reaching a particular target (red sphere) in a 2D square space. The target is placed at random at the start of each episode. So,  the action belongs to continuous space. Implementation of Deep Deterministic Policy Gradients (DDPG) algorithm from \cite{lillicrap2019continuous} is used to train the Actor-Critic agent. This policy $\pi \colon \mathcal{S} \to \mathcal{A}$ takes the state of a current environment and maps it to a predicted action. The Actor and Critic networks are made up of two hidden layers of fully-connected neural networks with $h_1 = 400$ and $h_2 = 300$ neurons with \textit{relu} activation functions. To match the constraints of available actions, the Actor's output is also scaled by \textit{tanh} activation multiplied by the maximum allowed action.

\begin{figure}
\centering
\begin{subfigure}{} 
	   \includegraphics[width=0.45\textwidth]{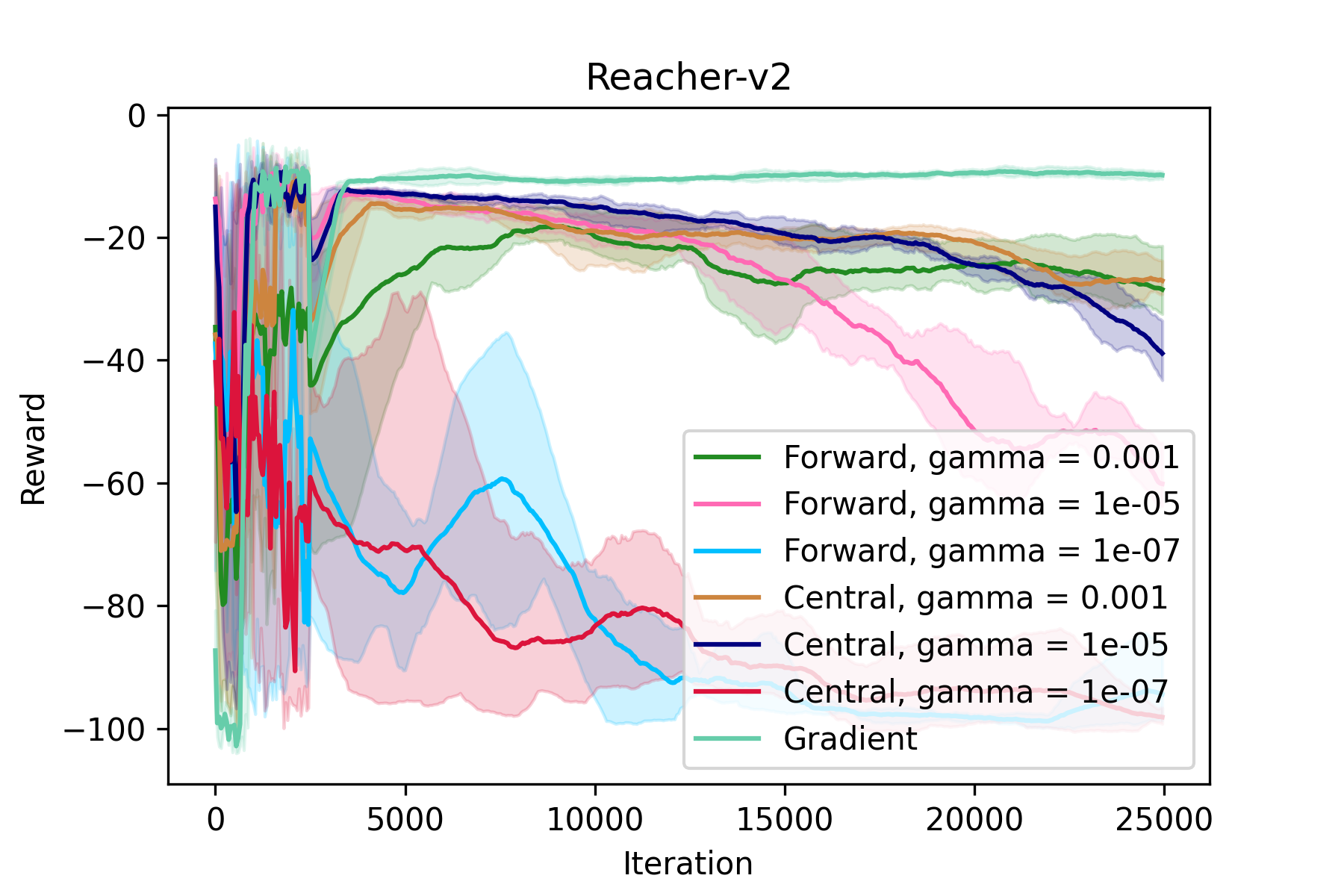}
\end{subfigure}
\hfill 
\begin{subfigure}{}
	   \includegraphics[width=0.45\textwidth]{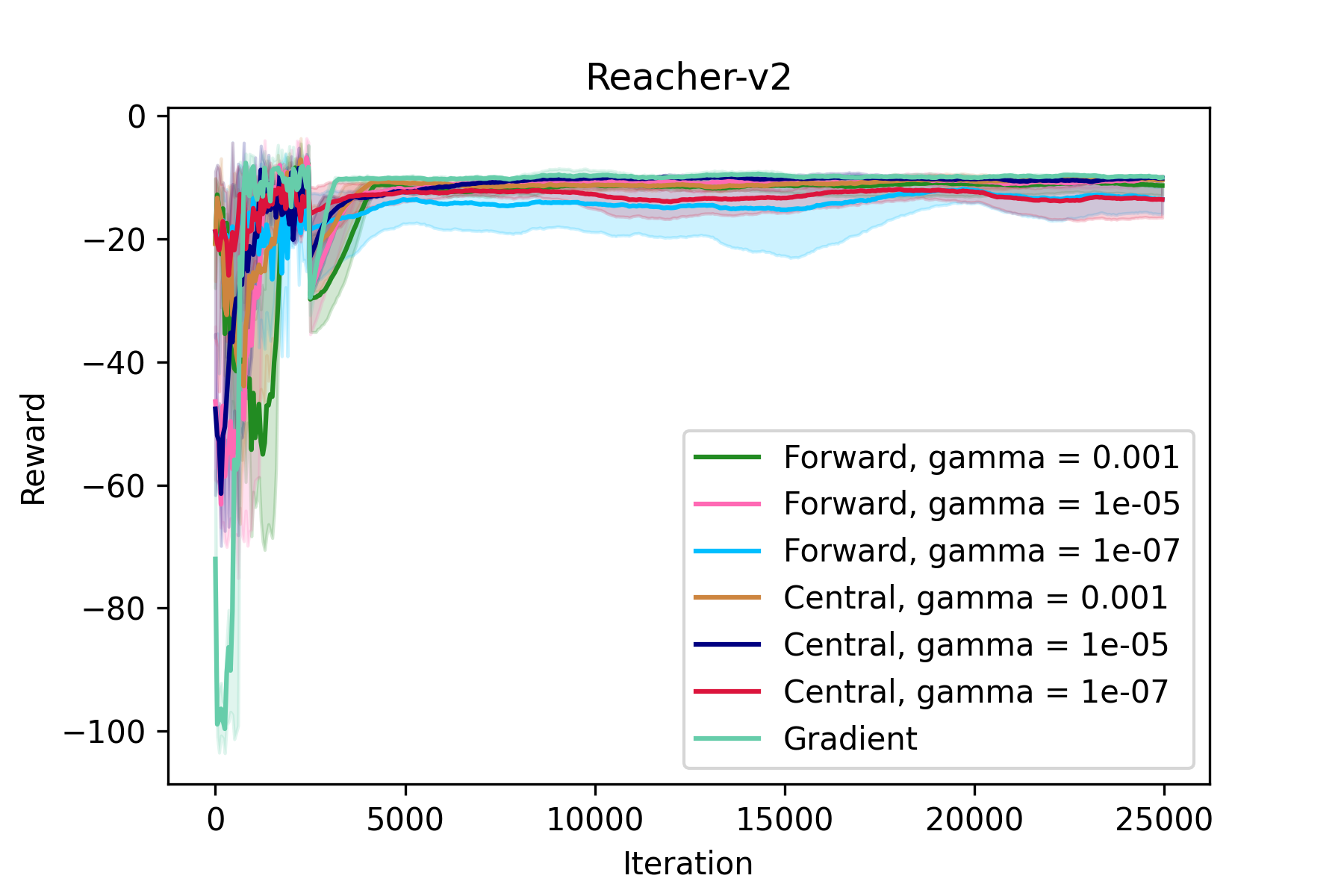} 
\end{subfigure}
\caption{Left: Actor's reward for ADAM with Forward and Central with various $\gamma$ and true gradient ADAM, where $lr = 0.001$.\\
Right: Actor's reward for ADAM with Forward and Central with various $\gamma$ and true gradient ADAM, where $lr = 0.0001$.}
\label{ap:fig_taus}
\end{figure}

\begin{figure}
\centering
\begin{subfigure}{}
	   \includegraphics[width=0.45\textwidth]{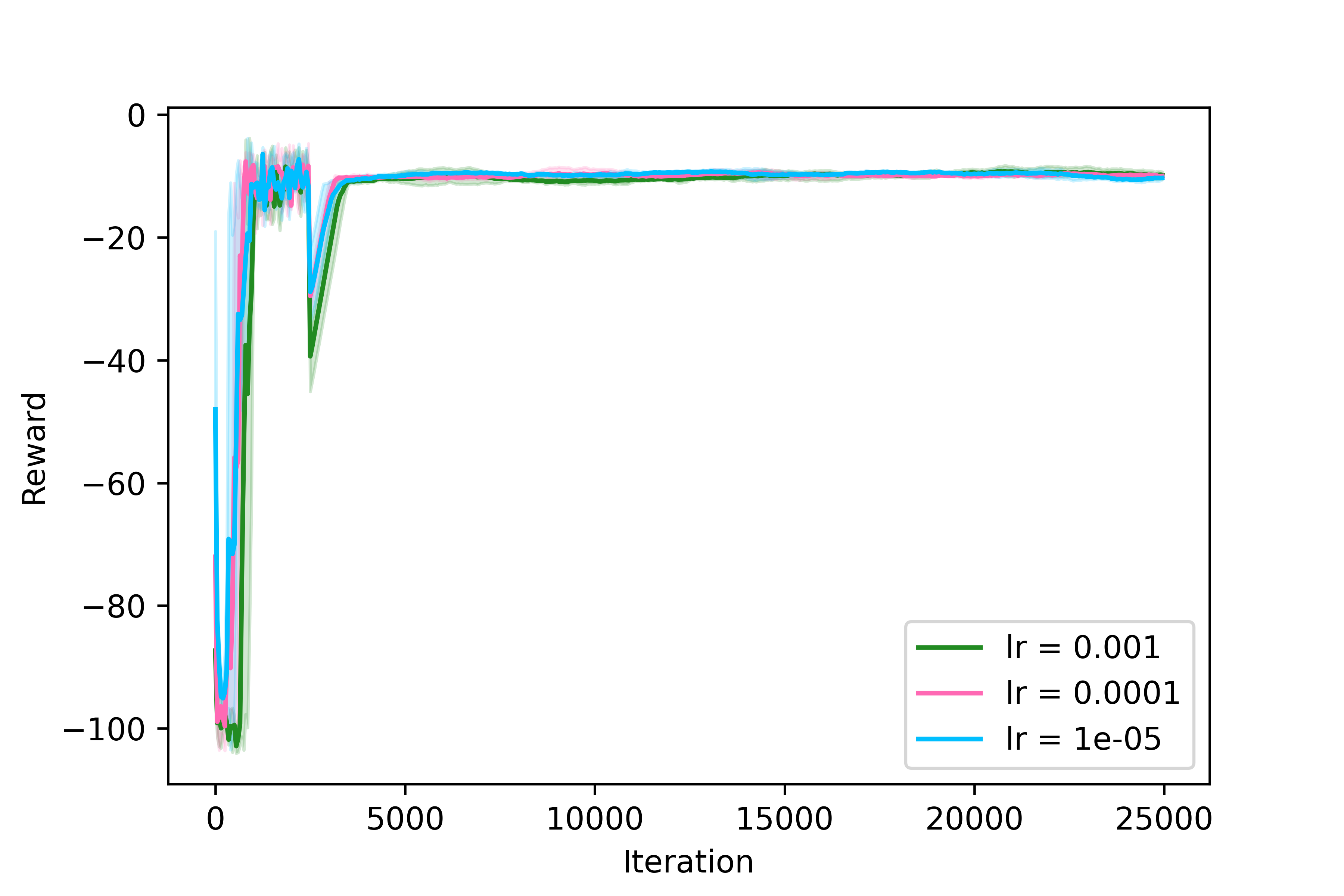} 
\end{subfigure}
\begin{subfigure}{}
	   \includegraphics[width=0.45\textwidth]{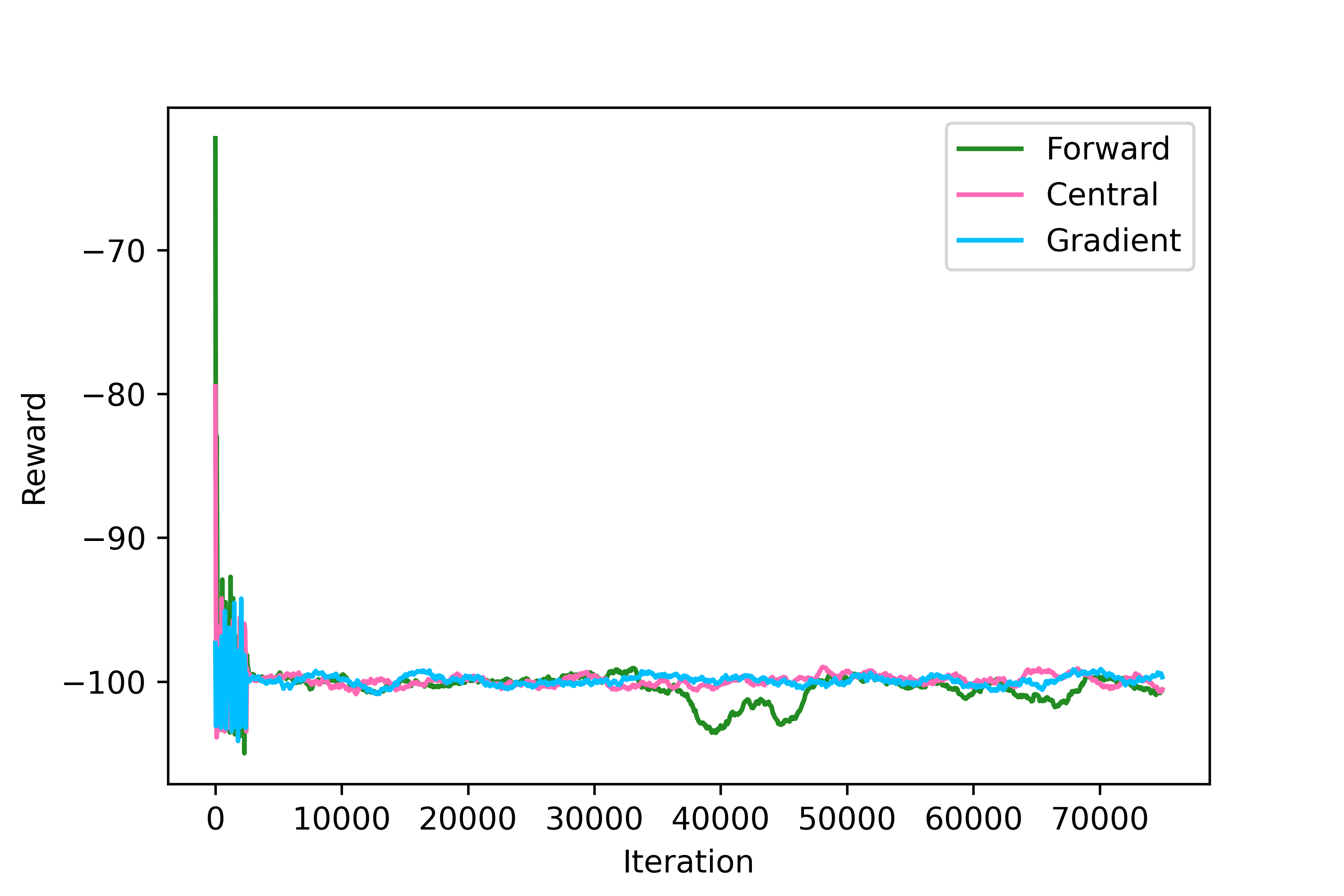} 
\end{subfigure}
\caption{Left: Actor's reward for true gradient ADAM with different learning rates.\\
Right: Actor's reward for $25000$ iterations of ADAM with Forward and Central finite difference with various $\gamma= 1e-5$ and true gradient ADAM, where $lr = 0.01$.}
\label{ap:fig_grads}
\end{figure}

\begin{figure}
\centering
\begin{subfigure}{}
	   \includegraphics[width=0.3\textwidth]{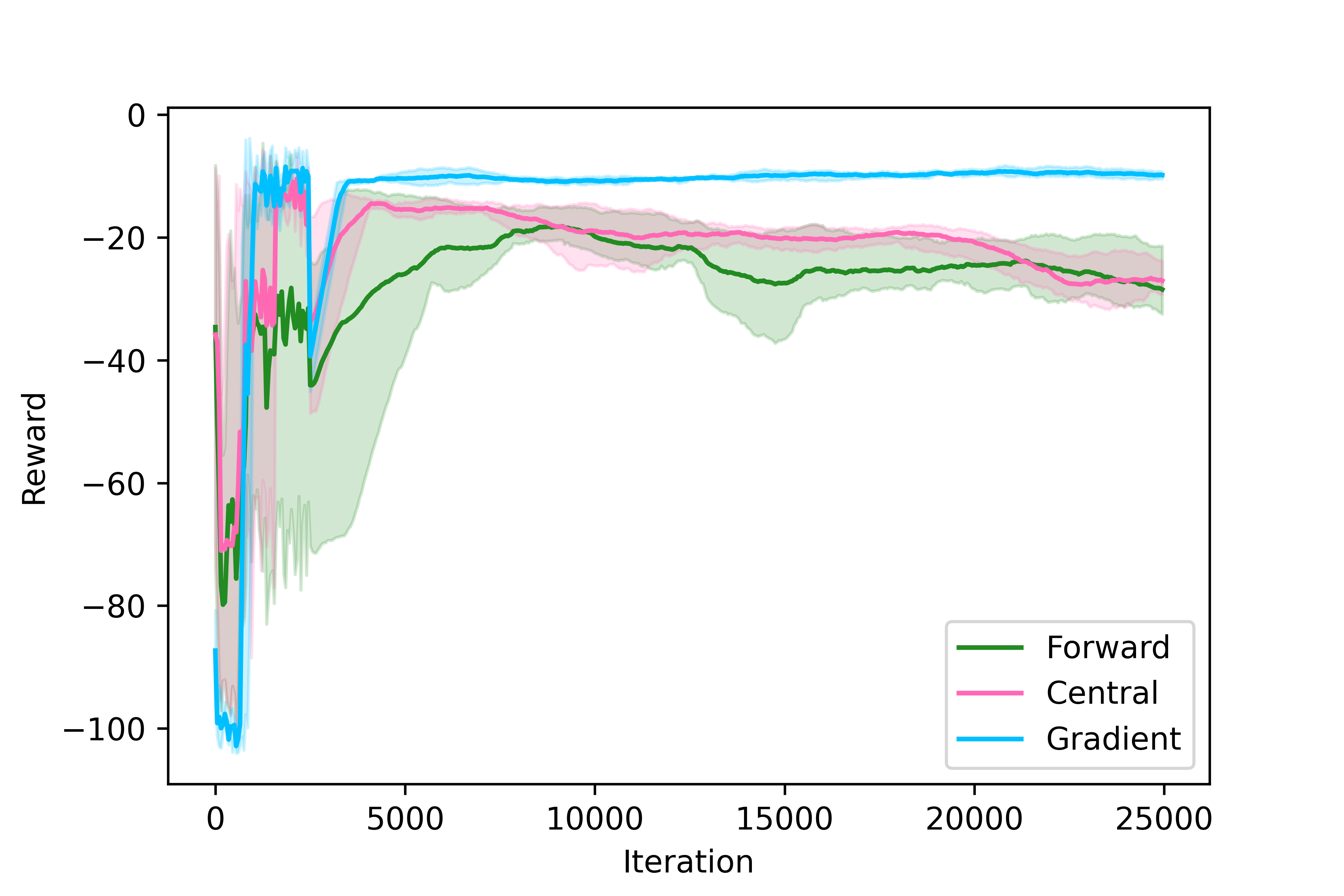}
\end{subfigure}
\begin{subfigure}{}
	     \includegraphics[width=0.3\textwidth]{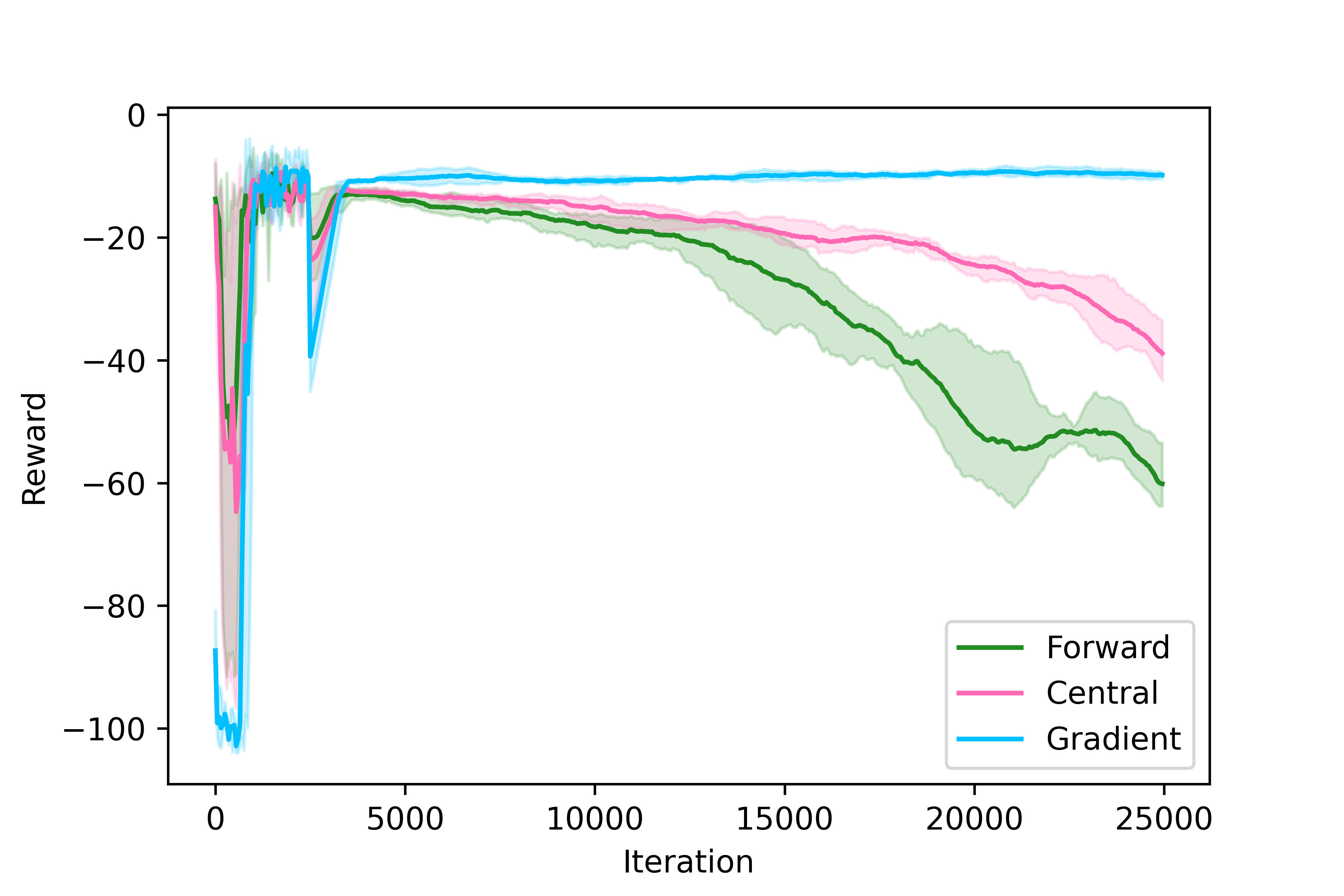}
\end{subfigure}
\begin{subfigure}{}
	    \includegraphics[width=0.3\textwidth]{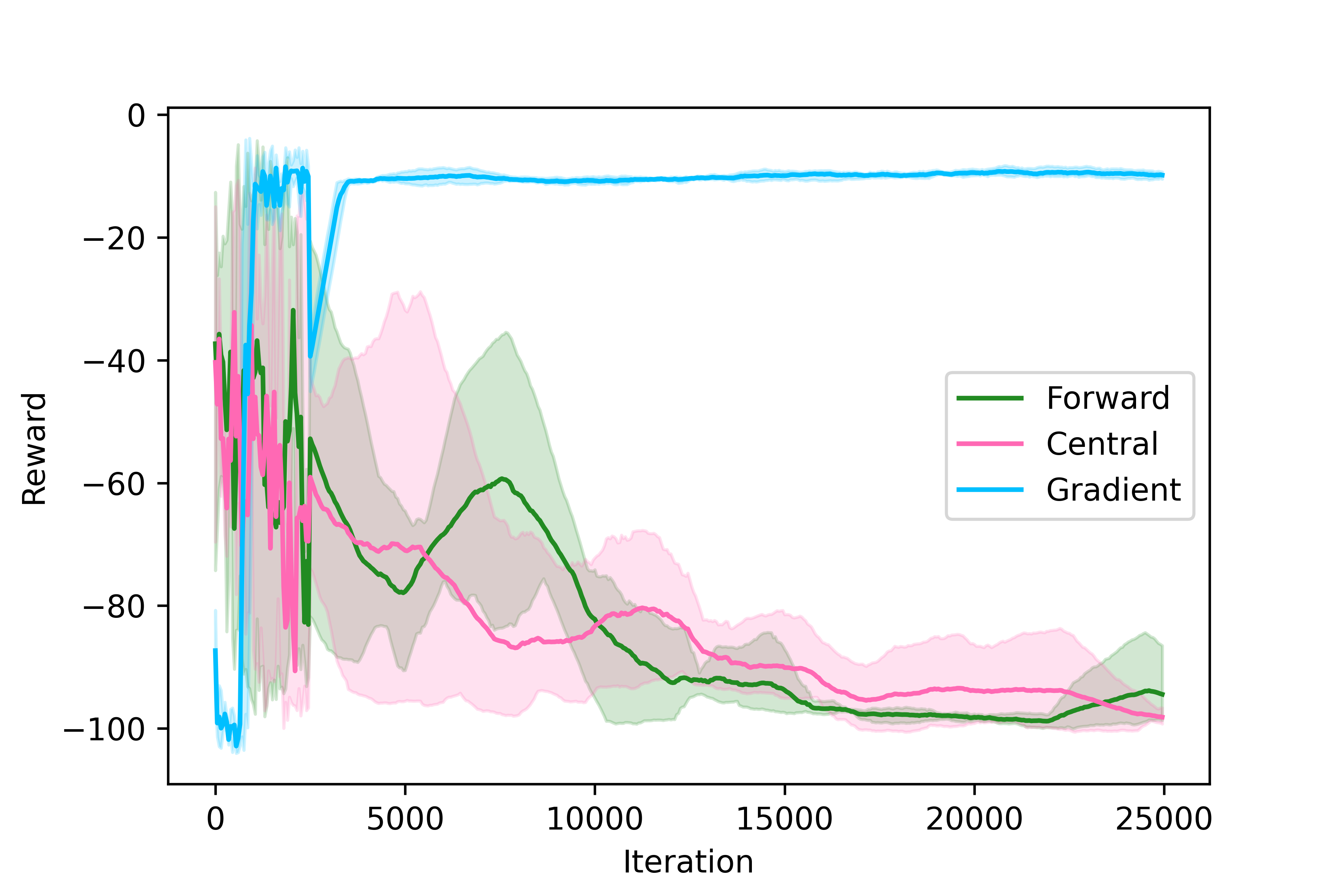}
\end{subfigure}
\caption{Left: Actor's reward for ADAM with Forward and Central for $\gamma= 1e-3$ and true gradient ADAM, where $lr = 0.001$.\\
Middle: Actor's reward for ADAM with Forward and Central for for $\gamma= 1e-5$ and true gradient ADAM, where $lr = 0.001$.\\
Right: Actor's reward for ADAM with Forward and Central for $\gamma= 1e-7$ and true gradient ADAM, where $lr = 0.001$.}
\label{ap:fig_lr001}
\end{figure}

Figure \ref{ap:fig_taus} shows the dependence of ZO methods on $\gamma$. One can see that: for a fitted learning rate $=0.0001$ at the right side of \ref{ap:fig_taus}, methods converge for all $\gamma$; for a small learning rate =1e-5 at \ref{rl_main}, methods converge for all $\gamma$ except the smallest $\gamma=$1e-7; for a big learning rate $=0.001$ at the left side of \ref{ap:fig_taus}, Gradient method converges but all ZO methods collect errors and slowly diverge; for a huge learning rate $=0.01$ at the right side of \ref{ap:fig_grads}, all methods diverge. As a result, we find a regime where ZO methods are stable and very close to gradient methods, but also, we find a regime where gradient method converges and ZO methods diverge. The left side of \ref{ap:fig_grads} shows that all gradient methods in our experiments converge to the same level. Figure \ref{ap:fig_lr001} is made for a detailed look on regime with a huge learning rate $=0.01$ at the left side of \ref{ap:fig_taus}.

\subsection{Additional Experiments: Robust Linear Regression}
\label{appendix_l1}
For the learning, we divide dataset "abalone scale" into two parts, where the part for training is $3500$ samples. The dimension of features equals to $8$. So, it is a very small dataset. We take dataset with such small dimension to compare Forward and Central with their variants of coordinate steps that depend on dimension. 

\textit{The Central Coordinate} is a zero-order finite difference such that we take  \eqref{sg} approximation $d$ times by each coordinate. Hence, we get $2d$ function computations per each step. 

\textit{The Forward Coordinate} is a zero-order finite difference such that we take  \eqref{forward} approximation $d$ times by each coordinate. Hence, we get $d+1$ function computations per each step.

As a result, coordinate steps are more accurate approximation of the gradient but also, they are more expansive computationally. We also add some graphics for Robust Linear Regression with different parameters.

\begin{figure}
\centering
\begin{subfigure}{}
\includegraphics[width=0.45\textwidth]{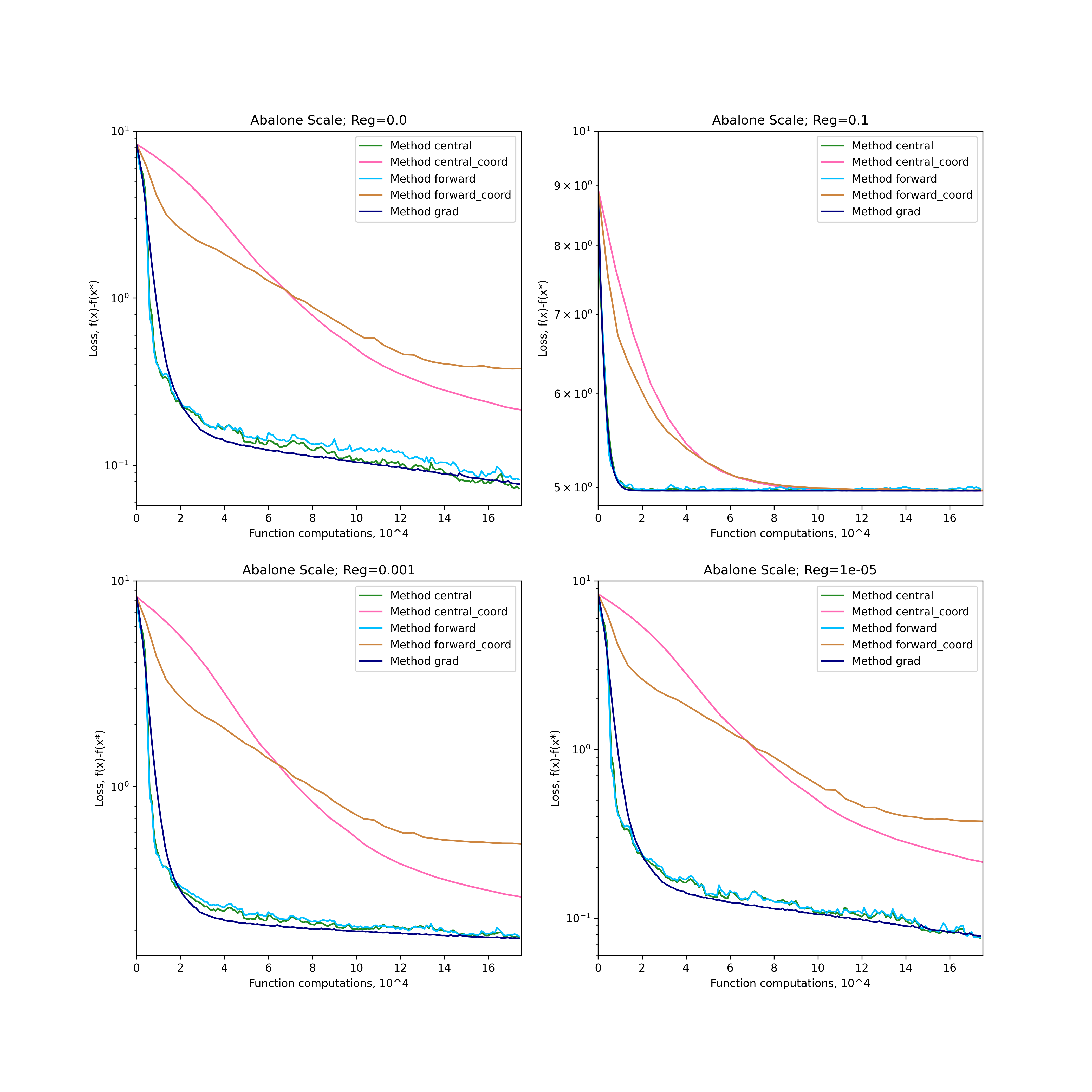}
\end{subfigure}
\begin{subfigure}{}
\includegraphics[width=0.45\textwidth]{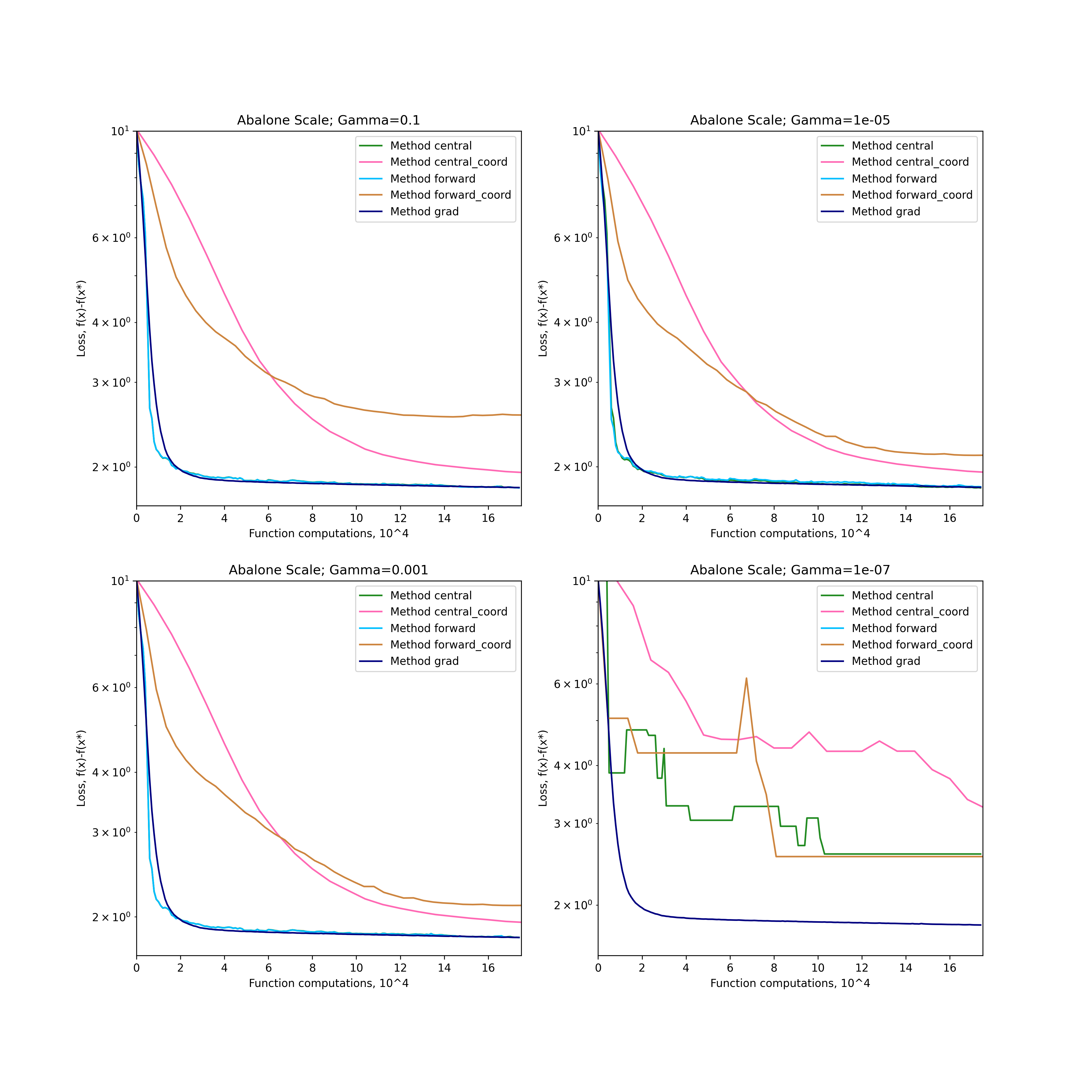}
\end{subfigure}

\caption{Left: Loss for \textit{abalone scale} dataset with $lr = 0.4$ \textit{batch size} $=100$, $\gamma=$1e-5, and different $\mu$.\\
Right: Loss for  \textit{abalone scale} dataset with $lr = 0.4$ \textit{batch size} $=100$, $\mu=0.$  and different $\gamma$.}
\label{ap:fig_l1_001}

\end{figure}

\begin{figure}
\centering
\begin{subfigure}{}
\includegraphics[width=0.45\textwidth]{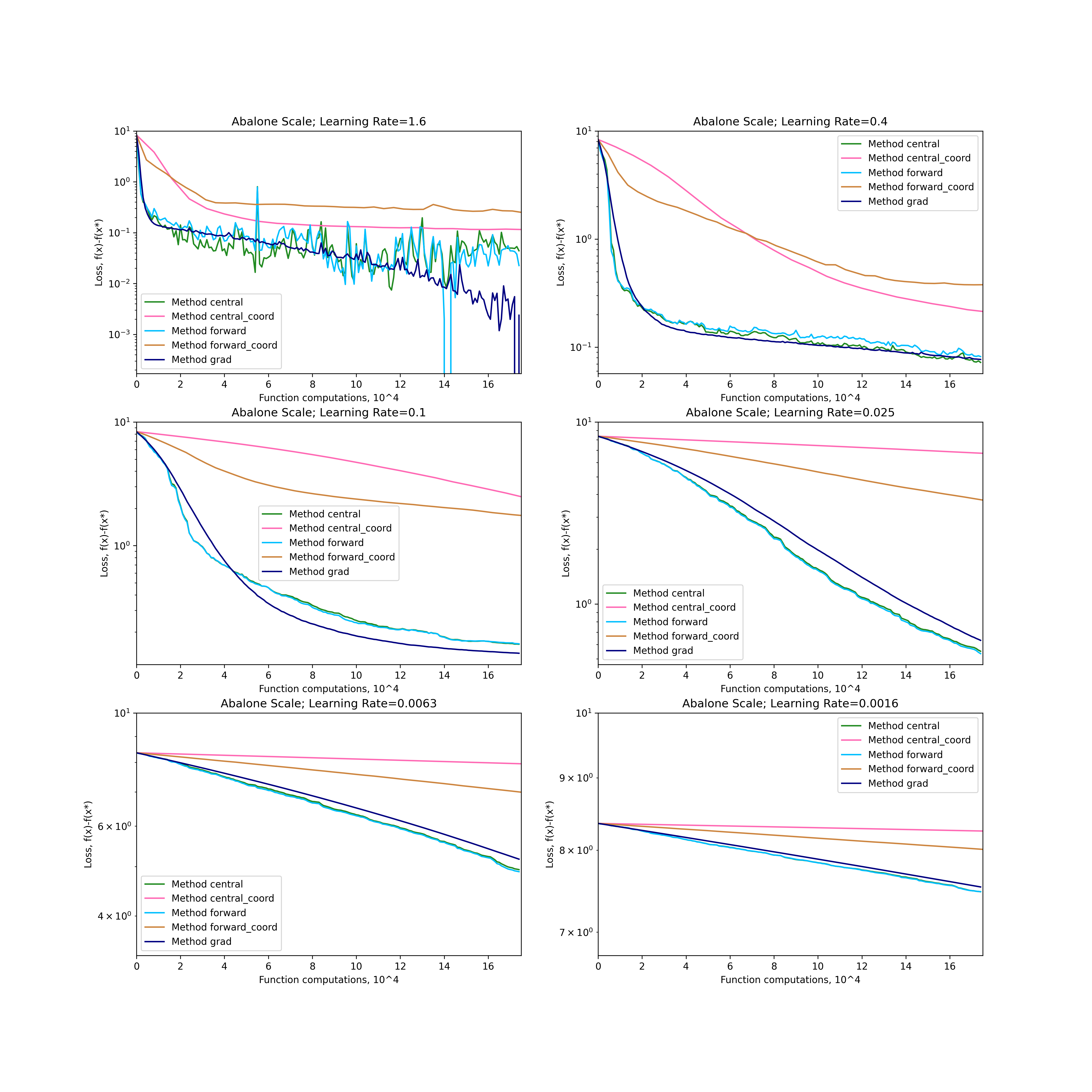}
\end{subfigure}
\begin{subfigure}{}
\includegraphics[width=0.45\textwidth]{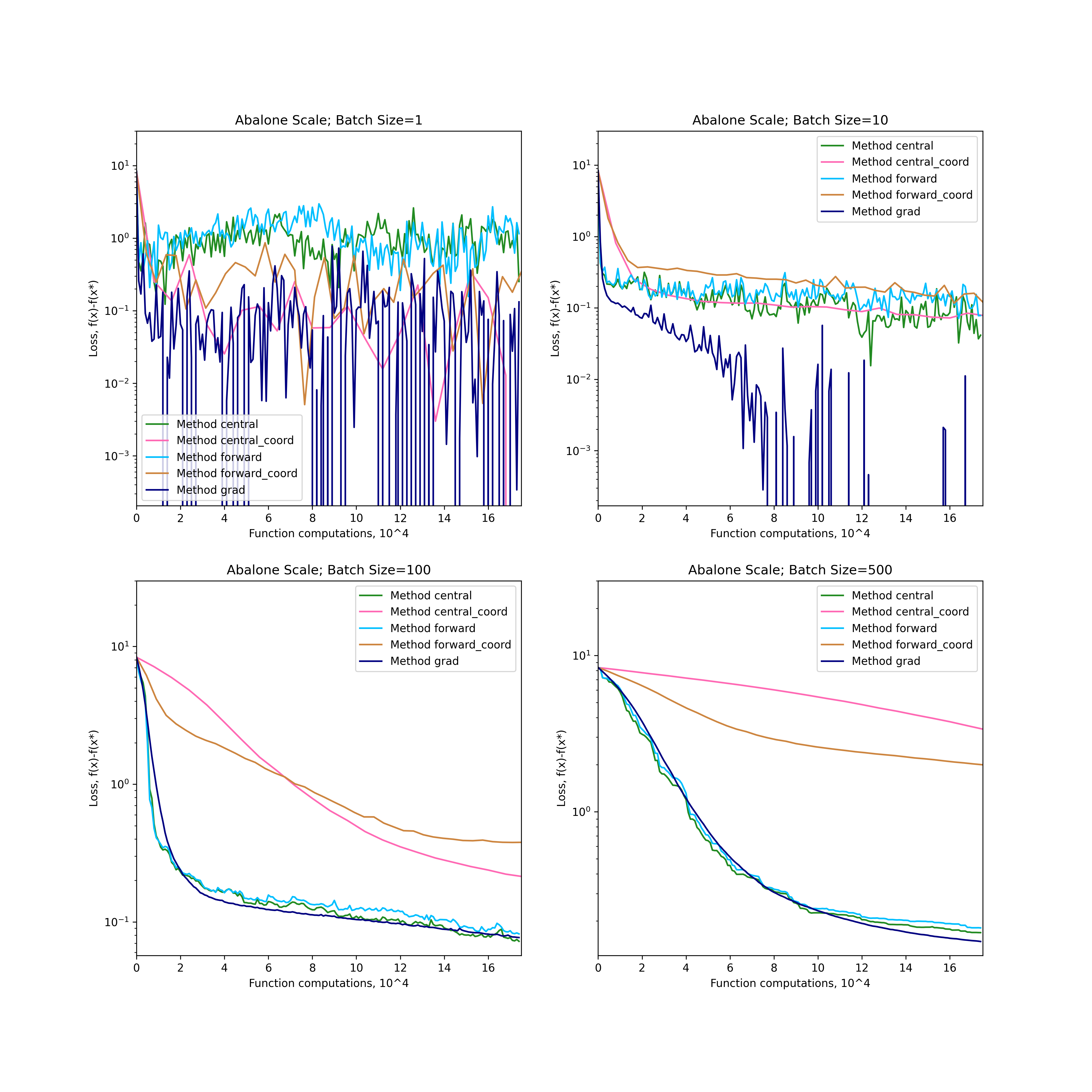}
\end{subfigure}
\caption{Left: Loss for \textit{abalone scale} dataset with $lr = 0.4$ \textit{batch size} $=100$, $\gamma=$1e-5, and different $lr$.\\
Right: Loss for  \textit{abalone scale} dataset with $lr = 0.4$, $\mu=0.$, $\gamma=$1e-5,  and different \textit{batch size}.}
\label{ap:fig_l1_002}

\end{figure}

\subsection{Additional Experiments: Support Vector Machine}
\label{appendix_svm}
In this subsection one can see additional graphics for SVM with different parameters. 

\begin{figure}
\begin{subfigure}{}
\includegraphics[width=0.45\textwidth]{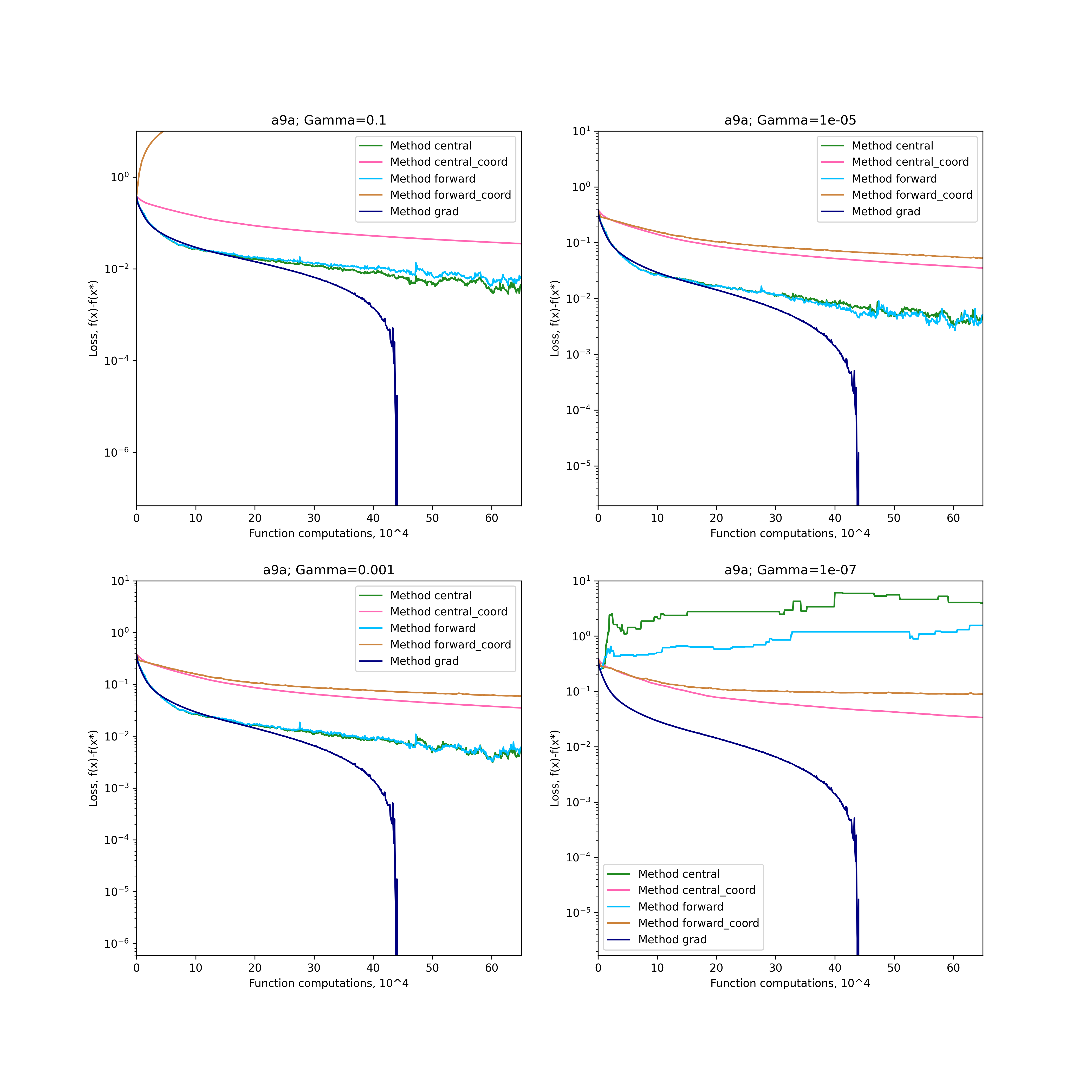}
\end{subfigure}
\begin{subfigure}{}
\includegraphics[width=0.45\textwidth]{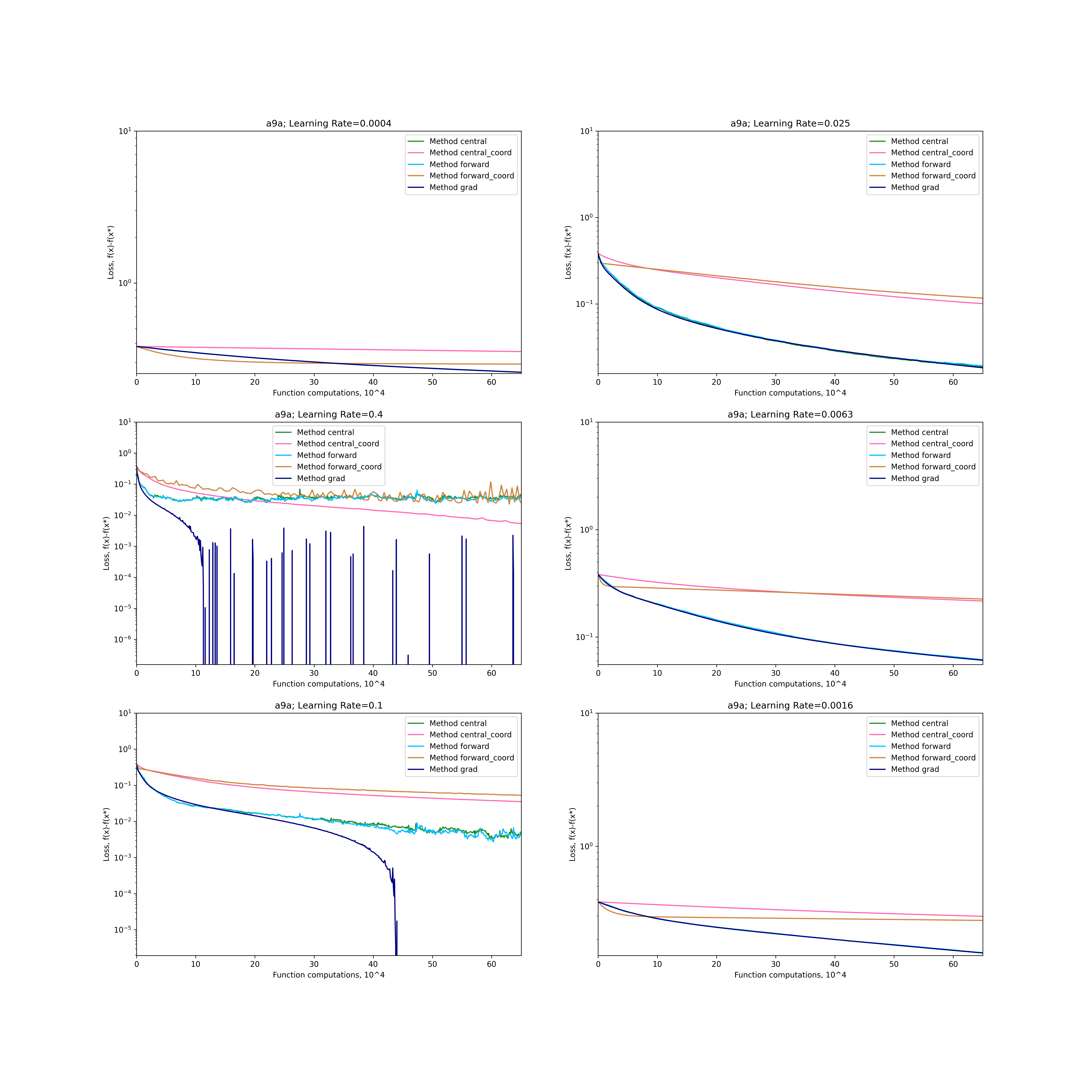}
\end{subfigure}
\caption{Left: Loss for \textit{a9a} dataset with $\mu=1e-5$, batch size $=100$, $lr = 0.1$ and different $\gamma$.\\
Right: Loss for \textit{a9a} dataset with $\mu=1e-5$, batch size  $=100$, $\gamma =$1e-5 and different $lr$.}

\label{ap:svm1}
\end{figure}

Figure \ref{ap:svm1} compares performance of ZO methods for $a9a$ dataset with various $\gamma$ (Left) and learninng rates (Right). It can be observed (Left) that \textit{Central} and \textit{Forward} methods converge faster than \textit{Central Coordinate} and \textit{Forward Coordinate} methods. Also, it is clear that both pairs of \textit{Coordinate} and \textit{Non-Coordinate} methods converge in a similar fashion. However, under extreme cases of $\gamma$ some methods do not converge at all. For instance, \textit{Forward Coordinate} method with $\gamma=0.1$, \textit{Central} and \textit{Forward} methods with $\gamma=1e^{-07}$. 
Comparison of ZO methods under different \textit{Learning rate} values (Right) also supports the above-mentioned corollary that in general \textit{Central} and \textit{Forward} methods converge faster than \textit{Central Coordinate} and \textit{Forward Coordinate} methods. However, it can be seen that larger values of learning rate introduces variance to ZO methods, although considerably smaller than for true gradient. It is also worth noting that small learning rate values result \textit{Central} and \textit{Forward} methods to converge with the same rate.

\begin{figure}
\centering
\begin{subfigure}{}
\includegraphics[width=0.85\textwidth]{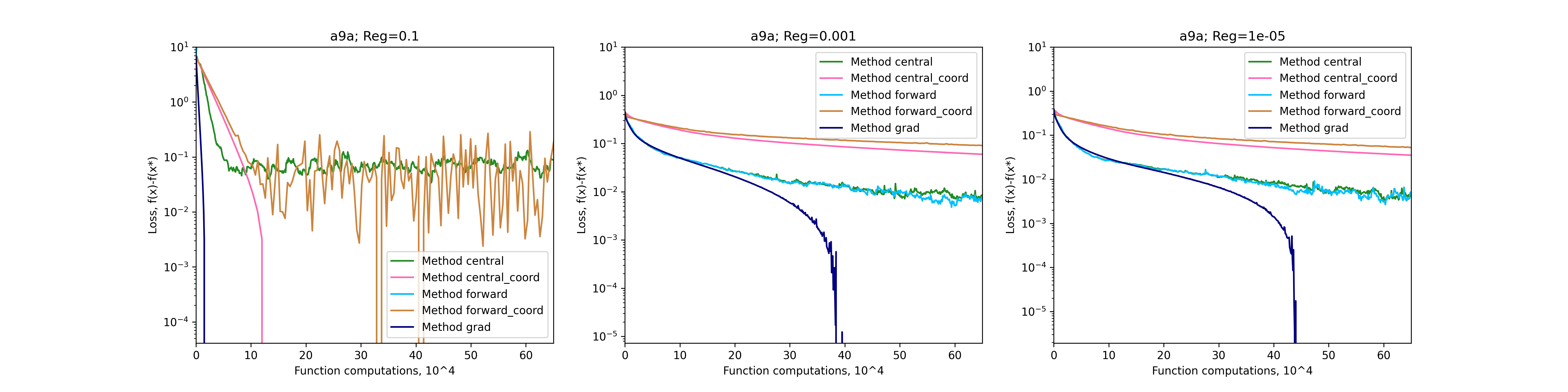}
\caption{Loss for \textit{a9a} dataset with $\mu=$1e-5, batch size $=100$, $lr = 0.1$, $\gamma =$1e-5 and different $\mu$.}
\end{subfigure}
\begin{subfigure}{}
\includegraphics[width=0.85\textwidth]{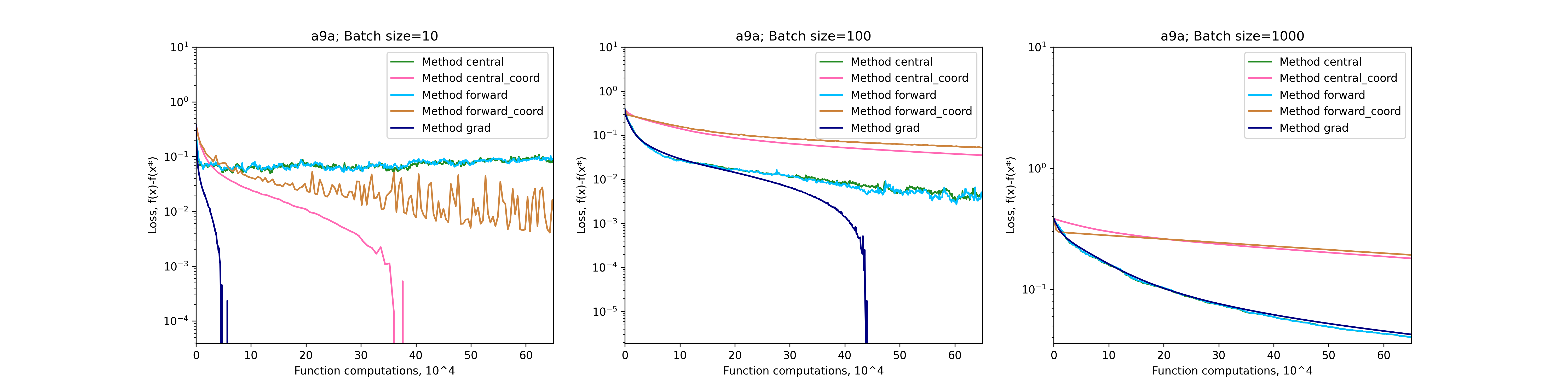}
\end{subfigure}
\caption{ Loss for \textit{a9a} dataset with $\mu=$1e-5, $lr = 0.1$, $\gamma =$1e-5 and different batch size.}
\label{ap:svm2}
\end{figure}

\subsection{Reducing Variance Under Batching In $p$-norm}
\label{appendix1}

\begin{lemma}
\label{lemma1}
If $\sigma$ is determined from
$$\mathbb{E}_{e \sim RS_2^d(1)} \left[ \exp \left( \frac{\left\| d\dfrac{f(x+ \gamma e) - f(x-\gamma e) }{2\gamma} e  - \nabla f_{\gamma}(x)\right\|_q^2}{\sigma^2} \right) \right] \leq \exp(1),$$
then the variance of batched gradient 
\begin{eqnarray*}
     \sigma_B^2 = \text{Var}\left[ \nabla^B f_{\gamma}(x, \{e_i\}_{i=1}^B) \right] \leq \dfrac{\sigma^2}{B} \cdot (2\chi(p, d) + \sqrt{3 \pi \chi(p, d)} + 3) 
     = \dfrac{\sigma^2}{B} \lambda(p,d),
\end{eqnarray*}
where $\chi(p, d) = \min{\left\{q-1, 2\ln{d} \right\}}$,  $\frac{1}{p}+\frac{1}{q}=1$.
\end{lemma}

\textbf{Proof:}

The batched gradient for the function $f_{\gamma}(x) = \mathbb{E}_u \left[ f(x + \gamma u) \right]$ ($u \in RB_2^d(1)$) is:
    \begin{eqnarray*}
    \label{f2}
        \nabla^B f_{\gamma}(x, \{e_i\}_{i=1}^B) =  \frac{1}{B} \sum\limits_{i=1}^B \left[ \dfrac{d}{2\gamma} \left( f(x+\gamma e_i) - f(x-\gamma e_i )\right)e_i \right],
    \end{eqnarray*} 
    where $e_i \in RS_2^d(1)$ are i.i.d.
    
    For simplicity we denote $s_i =\left[ \dfrac{d}{2\gamma} \left( f(x+\gamma e_i) - f(x-\gamma e_i )\right)e_i \right]$,  so
     $\nabla^B f_{\gamma}(x, \{e_i\}_{i=1}^B) = \dfrac{1}{B} \sum\limits_{i=1}^B s_i$, where $s_i$ are i.i.d.
     We denote $\overset{\circ}{s}_i = s_i - \mathbb{E} s_i$.

From Theorem $2.1$ from \cite{juditsky2008large}, we have
\begin{eqnarray*}
     \mathbb{P} \left( \left\| \sum\limits_{i=1}^B \overset{\circ}{s}_i \right\|_q \geq \left(\sqrt{2\chi} + \sqrt{2}\beta \right)\sqrt{B} \sigma \right) \leq e^{-\beta^2/3},
\end{eqnarray*}
where $\xi \leq  \min{\left( 2 \ln{d}, q-1 \right)}$ (example $3.2$ from \cite{juditsky2008large}), $q \geq 2$

\begin{eqnarray*}
     \text{Var} \left[ \nabla^B f_{\gamma}(x, \{e_i\}_{i=1}^B )\right] &=& \mathbb{E} \left[ \left\| \dfrac{1}{B} \sum\limits_{i=1}^B \overset{\circ}{s}_i \right\|_q^2 \right] = \int\limits_{t=0}^{+\infty} \mathbb{P} \left( \left\| \dfrac{1}{B} \sum\limits_{i=1}^B \overset{\circ}{s}_i \right\|_q  \geq \sqrt{t} \right)dt \nonumber\\
     &&= \int\limits_{t=0}^{2\chi B \sigma^2} \mathbb{P} \left( \left\| \dfrac{1}{B} \sum\limits_{i=1}^B \overset{\circ}{s}_i \right\|_q  \geq \sqrt{t} \right)dt + \int\limits_{t = 2\chi B \sigma^2}^{+\infty} \mathbb{P} \left( \left\| \dfrac{1}{B} \sum\limits_{i=1}^B \overset{\circ}{s}_i \right\|_q  \geq \sqrt{t} \right)dt. \end{eqnarray*}

     Substituting  $\sqrt{t} = \left( \sqrt{2\chi} + \sqrt{2}\beta \right)\sqrt{B}\sigma$, we obtain:
     \begin{eqnarray*}
          \text{Var} \left[ \nabla^B f_{\gamma}(x, \{e_i\}_{i=1}^B )\right] &\leq& 2\chi B \sigma^2 + \int\limits_{\beta = 0}^{+\infty} \left[\mathbb{P} \left( \left\| \dfrac{1}{B} \sum\limits_{i=1}^B \overset{\circ}{s}_i \right\|_q \geq \left( \sqrt{2\chi} + \sqrt{2}\beta \right)\sqrt{B} \sigma \right)- 2B \sigma^2(\sqrt{\chi} + \beta)\right]d\beta \nonumber\\
          &&\leq 2\chi B \sigma^2 + 2B\sigma^2 \int\limits_{\beta = 0}^{+\infty} e^{-\beta^2/3} (\sqrt{\chi} + \beta)d\beta = 2B\sigma^2 \left(  \chi + \dfrac{\sqrt{3\pi}}{2}\sqrt{\chi} + \dfrac{3}{2} \right).
     \end{eqnarray*}
     
     Finally, 
     \begin{eqnarray*}
          \sigma_B^2 = \text{Var} \left[ \nabla^B f_{\gamma}(x, \{e_i\}_{i=1}^B )\right] = \text{Var} \left[ \dfrac{\sum\limits_{i=1}^B s_i}{B} \right] \leq \dfrac{\sigma^2}{B} \cdot (2\chi(q,d) + \sqrt{3\pi\chi(q,d)} + 3) = \dfrac{\sigma^2}{B} \lambda(p,d),
     \end{eqnarray*}
     where 
     \begin{equation}\label{lambda_def}
         \lambda(p,d) = 2\chi(p,d) + \sqrt{3\pi\chi(p,d)} + 3, \;
      \chi(p,d) = \min{\left( q-1, 2\ln{d} \right)}.
    \end{equation}
    
     For example, for $p = 1$ ($q=\infty$) we have $\sigma_B^2 \leq \dfrac{\sigma^2}{B} \cdot O \left( \ln{d} \right)$ and for $p=2$ we have $\sigma_B^2 \leq \dfrac{9\sigma^2}{B}$, which coincides with the well known property $\sigma_B^2 = \dfrac{\sigma^2}{B}$ up to a numerical constant.


\subsection{Proof Of \textbf{Theorem \ref{main_properties} (Properties Of $f_{\gamma}$})}

\label{main_properties_proof} 

For all $x,y\in \ag{Q}$ 

\begin{itemize}
    \item[$1.$] $f(x) \leq f_{\gamma}(x) \leq f(x) + \gamma M_{\ag{2}}$;
    
    \item[$2.$] $f_{\gamma}(x)$ is $M$-Lipschitz: $$|f_{\gamma}(y) - f_{\gamma}(x) | \leq M \| y - x\|_p;$$
    
    \item[$3.$] $f_{\gamma}(x)$ has $L = \dfrac{2\sqrt{d} M}{\gamma}$-Lipschitz gradient:
    $$\|\nabla f_{\gamma}(y) - \nabla f_{\gamma}(x) \|_q \leq L \| y - x\|_p.$$ 
\end{itemize}
where $q$ us such that $1/p + 1/q = 1$. 

\textbf{Proof:}

For the first point, we have:

For the first inequality, we use the convexity of function $f(x)$
\begin{eqnarray*}
    f_{\gamma}(x) = \mathbb{E} \left[ f(x + \gamma u) \right]  \geq \mathbb{E} \left[ f(x) + \left\langle \ag{\nabla f}(x), \gamma u \right\rangle \right] = \mathbb{E}\left[f(x) \right]  = f(x)
\end{eqnarray*}

For the second inequality:
\begin{eqnarray*}
    \left| f_{\gamma}(x) - f(x) \right| = \left| \mathbb{E} \left[ f(x+\gamma u) \right] - f(x)\right| \leq \mathbb{E} \left[ \left| f(x + \gamma u) - f(x)   \right| \right] \leq  \an{\mathbb{E} \left[ M_2 \cdot \| \gamma u \|_2\right]} \leq \gamma M_{\an{2}},
\end{eqnarray*}
\an{using the fact that $f$ is $M_2$-Lipshcitz}.

For the second point:
\begin{eqnarray*}
    | f_{\gamma}(y) - f_{\gamma}(x) | = |\mathbb{E} \left[ f(y + \gamma u) - f(x + \gamma u) \right] | \leq \mathbb{E} | f(y + \gamma u)  - f(x + \gamma u)| \leq M \|y - x \|_p.
\end{eqnarray*}

In the third point, applying Lemma $11$ from \cite{duchi2012randomized},  we have:
\begin{eqnarray*}
    \left\| \nabla f_{\gamma}(y) - \nabla f_{\gamma}(x) \right\|_q &=& \left\| \nabla \mathbb{E}_{Z \sim B^d_2(\gamma)} \left[ f(y + Z)  \right] - \nabla \mathbb{E}_{Z \sim B^d_2(\gamma)} \left[ f(x + Z)  \right] \right\|_q  \nonumber\\ &&= \left\|  \mathbb{E}_{Z \sim B^d_2(\gamma)} \left[ \nabla f(y + Z)  \right] -  \mathbb{E}_{Z \sim B^d_2(\gamma)} \left[ \nabla f(x + Z)  \right] \right\|_q \nonumber\\ &&\leq M \underbrace{\int | \mu(z-y) - \mu(z-x) |dz }_{I_1},
\end{eqnarray*}
where  $\mu(x) = \dfrac{1}{V(B^d_2(\gamma))} \cdot \mathbb{I} \left( x \in B^d_2(\gamma) \right)$. Note that $f(x)$ is not assumed to be differentiable but the Lebesgue measure of the set where the convex function is not differentiable is equal to zero.

Using the bound for Integral $I_1$ from Lemma $8$ from \cite{YOUSEFIAN201256} and the fact, that $$\lim\limits_{d \rightarrow \infty} \dfrac{\kappa \dfrac{d!!}{(d-1)!!}}{\sqrt{d}} = \dfrac{\sqrt{\pi}}{2},$$ we obtain 
\begin{eqnarray*}
    \left\| \nabla f_{\gamma}(y) - \nabla f_{\gamma}(x)  \right\|_q &\leq& \dfrac{\sqrt{d} M}{\gamma} \sqrt{\dfrac{2}{\pi}} \left\| y-x \right\|_2. 
\end{eqnarray*}

Since $\|y - x \|_2 \leq \|y - x \|_p$ for $p \in [1, 2]$ and $\pi > 2$, we obtain:
\begin{eqnarray*}
    \left\| \nabla f_{\gamma}(y) - \nabla f_{\gamma}(x)  \right\|_q &\leq& \dfrac{\sqrt{d} M}{\gamma}  \left\| y-x \right\|_p.
\end{eqnarray*}

\subsection{Proof Of \textbf{Theorem \ref{main_properties2} (Properties Of $\nabla f_{\gamma}(x,e)$})}
\label{theorem_22_proof}
For all $x\in Q$

 \begin{itemize} 
    \item Unbiased:
    $\mathbb{E}_e \left[\nabla f_{\gamma}(x, e) \right] = \nabla f_{\gamma}(x)$;
    
    \item  Bounded variance (second moment):
    $$ \mathbb{E}_e \left[ \| \nabla f_{\gamma}(x, e) \|^2_q \right] =
    \kappa(p,d)\cdot\left(d M_2^2 + \dfrac{d^2\Delta^2}{\gamma^2} \right),$$
   
   where $1/p + 1/q = 1$ and
   \begin{equation*}
\kappa(p,d) = O\left(\sqrt{\mathbb{E}_e\|e\|_q^4}\right) = 
 \begin{cases}
   O(1), \;  p = 2\\
   O\left((\ln d)/d\right), \; p = 1.
 \end{cases}
\end{equation*}
\end{itemize}
If $\Delta$ is sufficiently small, then 
 $$\mathbb{E}_e \left[ \| \nabla f_{\gamma}(x, e) \|^2_q \right]\lesssim 2\kappa(p,d)dM_2^2.$$

\textbf{Proof:}

For the first point, substitute $z = \gamma u$, then, according the definition of $f_{\gamma}(x)$
\begin{eqnarray*}
    f_{\gamma}(x) = \dfrac{1}{V\left( B_2^d (\gamma) \right)} \int\limits_{\| z \|_2 \leq \gamma} f(x+z)dz.
\end{eqnarray*}

Since $f(x)$ is continuous, $f_{\gamma}(x)$ is continuously differentiable and its gradient can be found from \cite{nemirovsky1983problem} (see  formula 3.2 in chapter \textbf{9.3.2}):
\begin{eqnarray*}
    \nabla f_{\gamma}(x) = \dfrac{1}{V\left( B_2^d (\gamma) \right)} \int\limits_{\| z \|_2 = \gamma} f(x+z) \dfrac{z}{\| z\|_2}dS_{\gamma}(z),
\end{eqnarray*}
where $dS_{\gamma}(z)$ is an element of a spherical surface of radius $\gamma$

After normalization to the normalized area (the area of the whole sphere is taken $1$) we have integration with respect to  a uniformly distributed probability $d \sigma(e)$ on $S_1$
\begin{eqnarray*}
    \nabla f_{\gamma}(x) = \dfrac{d}{\gamma} \int\limits_{\| e \|_2 = 1} f(x + \gamma e) d \sigma(e) = \mathbb{E}_{e \sim RS_1^d(0)} \left[ \dfrac{d f(x + \gamma e) \cdot e}{\gamma} \right].
\end{eqnarray*}

Since $f(x + \gamma e) \cdot e$ has the same distribution as $f(x - \gamma e) \cdot e$, we also get:

\begin{eqnarray*}
    \nabla f_{\gamma}(x) = \mathbb{E}_{e \sim RS_1^d(0)} \left[ \dfrac{d \left(f(x + \gamma e) - f(x - \gamma e) \right) }{2\gamma} e \right] = \mathbb{E}_e \left[ \nabla f_{\gamma} (x, e) \right].
\end{eqnarray*}

The second point is proved in Lemma 2 from \cite{beznosikov2020derivative} (the second statement).

\subsection{Proof Of \textbf{Theorem \ref{FGM}}}
\label{th2.4proof}

It is proven in \cite{lan2012optimal}, that algorithm after $N$ iteration gives accuracy for $f_{\gamma}(x):$
\begin{eqnarray*}
    \mathbb{E} \left[ f_{\gamma} (x_{ag}^{N+1}) - f(x_*(\gamma)) \right] \leq \dfrac{4L_{f_{\gamma}} R^2}{N^2} + \dfrac{4\sigma_B R }{\sqrt{N}},
\end{eqnarray*}
where 
$x_*(\gamma) = \underset{x \in Q_{\gamma}}{\operatorname{argmin}} f_{\gamma}(x)$, $L_{f_{\gamma}} = \dfrac{\sqrt{d}\an{\sqrt{M_2 M}}}{\gamma}$.

If we have $\dfrac{\varepsilon}{2}$-accuracy for the function $f_{\gamma}(x) $ with $\gamma = \dfrac{\varepsilon}{2M_{\an{2}}}$, then we have $\varepsilon$-accuracy for the function $f(x)$:
\begin{eqnarray*}
    f(x_{N+1}^{ag}) - f(x_*) \leq f(x_{N+1}^{ag}) - f(x_*(\gamma))  \leq f_{\gamma} (x_{N+1}^{ag}) + \gamma M_{\an{2}} - f_{\gamma}(x_*(\gamma)) \leq \dfrac{\varepsilon}{2} + \dfrac{\varepsilon}{2} = \varepsilon.
\end{eqnarray*}

To have $\dfrac{\varepsilon}{2}$-accuracy for $f_{\gamma}(x)$ we need 
\begin{eqnarray*}
    \dfrac{4L_{f_{\gamma}} R^2}{N^2} \leq \dfrac{\varepsilon}{4}
\end{eqnarray*}
and
\begin{eqnarray*}
    \dfrac{4\sigma_B R }{\sqrt{N}} \leq \dfrac{\varepsilon}{4}.
\end{eqnarray*}

Substituting $L_{f_{\gamma}}$ from Theorem \ref{main_properties} and $\sigma_B^2$ from Lemma \ref{lemma1}  we obtain:
\begin{eqnarray*}
    N = \dfrac{4\sqrt{2}  \an{\sqrt{M_2 M}} R}{\varepsilon} = O \left( \dfrac{ d^{1/4} \an{\sqrt{M_2 M}} R}{\varepsilon} \right),
\end{eqnarray*}
 and \begin{eqnarray*}
    B = \max{\left( 1, \dfrac{256 \lambda(p,d)\sigma^2 R^2}{\varepsilon^2 N} \right)} = \max{\left(1,  64\sqrt{2} \dfrac{\kappa(p, d) \lambda(p,d) d^{3/4} M_2^2 R}{M\varepsilon} \right)},
\end{eqnarray*}
 see $\lambda(p,d)$ in \eqref{lambda_def}. 
 
We obtain total number of oracle calls $T$: 
\begin{eqnarray*}
    && T = N \cdot B = \max{\left\{ \dfrac{ d^{1/4} \an{\sqrt{M_2 M}} R}{\varepsilon}, \dfrac{256 \lambda(p,d) \sigma^2 R^2} {\varepsilon^2} \right\}} =  \max{\left\{ 4\sqrt{2}\dfrac{ d^{1/4} \an{\sqrt{M_2 M}} R}{\varepsilon}, \dfrac{512\kappa(p, d) \lambda(p,d) d M_2^2 R^2} {\varepsilon^2} \right\}} \nonumber\\
    &&= \tilde{O} \left( \max{\left\{ \dfrac{ d^{1/4} \an{\sqrt{M_2 M}} R}{\varepsilon}, \dfrac{\kappa(p, d)  d M_2^2 R^2}{\varepsilon^2} \right\}} \right),
\end{eqnarray*}
as $\lambda(p,d) = O(\log(d)) = \tilde O(1)$, see \eqref{lambda_def}.

\subsection{$L_{xy}$ Estimate}
\label{L_xy_proof}

Applying Lemma $11$ from \cite{duchi2012randomized}, we have:
\begin{eqnarray*}
    \left\| \nabla_x f_{\gamma}(x, y_2) -  \nabla_x f_{\gamma}(x, y_1)\right\|_q &=& \left\| \nabla_x \mathbb{E}_{Z_x, Z_y}  \left[ f(x + Z_x, y_2 + Z_y)\right] -  \nabla_x \mathbb{E}_{Z_x, Z_y}  \left[ f(x + Z_x, y_1 + Z_y)\right]\right\|_q \nonumber\\
    &&= \left\| \mathbb{E}_{Z_x, Z_y} \left[ \nabla_x f(x + Z_x, y_2 + Z_y) - \nabla_x f(x + Z_x, y_1 + Z_y)  \right] \right\|_q \nonumber\\
    &&\leq M_x \underbrace{\int | \mu_y(z_y-y_2) - \mu(z_y-y_1) |dz }_{I_1(y)},
\end{eqnarray*}
where $\mu_y(y) = \dfrac{1}{V(B^{d_y}_2(\gamma_y))} \cdot \mathbb{I} \left( y \in B^{d_y}_2(\gamma_y) \right)$. Note that $f(x, y)$ is not assumed to be differentiable but the Lebesgue measure of the set where the convex-concave function is not differentiable is equal to zero.

Using the bound for Integral $I_1(y)$ from Lemma $8$ \cite{YOUSEFIAN201256} and the fact, that $$\lim\limits_{d \rightarrow \infty} \dfrac{\kappa \dfrac{d!!}{(d-1)!!}}{\sqrt{d}} = \dfrac{\sqrt{\pi}}{2},$$ we obtain:
\begin{eqnarray*}
     \left\| \nabla_x f_{\gamma}(x, y_2) -  \nabla_x f_{\gamma}(x, y_1)\right\|_q \leq \dfrac{M_x \sqrt{d_y}}{\gamma_y} \|y_2 - y_1 \|_2 \leq \dfrac{M_x \sqrt{d_y}}{\gamma_y} \|y_2 - y_1 \|_p,
\end{eqnarray*}
where $p \in [1, 2]$.

\subsection{Proof Of Theorem \ref{SFGM}}
\label{th3.1.proof}

Based on batched Accelerated gradient method Smoothing scheme gives gradient-free method with $$\tilde{O}\left(\frac{d^{1/4} \an{\sqrt{M_2 M}}}{\sqrt{\mu\varepsilon}}\right)$$ successive iterations and
$$\tilde{O}\left(\frac{\kappa(p,d)d M_2^2}{\mu\varepsilon}\right)$$
oracle calls, where $\kappa(p,d)$ defined in Theorem~\ref{FGM}. 

This result will be true for stochastic problem \eqref{stoch_problem} if $M_2$ is defined as $\mathbb{E}_{\xi}\|\nabla_x f(x,\xi)\|_2^2\le M^2_2$ for all $x\in Q_{\gamma}$.

\textbf{Proof:}
Below we use \textit{restarts scheme} described in \cite{juditsky2014deterministic}. For simplicity, we will denote $R$, $R_k$ distance from starting, current point to the solution in $p$-norm up to a $O(\ln d)$-factor in worth case.

Theorem \ref{FGM} proves that to achieve the error $\varepsilon_k$, we need $$N_{\varepsilon_k
}=\dfrac{4\sqrt{2} d^{1/4} \an{\sqrt{M_2 M }} R_k}{\varepsilon_k}$$ iterations and 
$$T(\varepsilon_k) = \max{\left\{ 4\sqrt{2}\dfrac{ d^{1/4} \an{\sqrt{M_2 M}} R_k}{\varepsilon_k}, \dfrac{512\kappa(p, d) \lambda(p,d) d M_2^2 R_k^2} {\varepsilon_k^2} \right\}}$$
oracle calls, see $\lambda(p,d)$ in \eqref{lambda_def}. 
We can take $$R_k = \sqrt{\dfrac{2\varepsilon_k}{\mu}}$$ as the function $f$ and consequently $f_{\gamma}$ are $\mu$-strongly convex. 
We take 
\begin{center}
$\varepsilon_k = \dfrac{\mu R^2}{2} \cdot 4^{-k}$, $\varepsilon_K = \varepsilon \; \implies \; K = \dfrac{\ln{\left( \dfrac{\mu R^2}{2\varepsilon}\right)}}{\ln{4}}$.
\end{center}
We obtain the number of iterations and the $k^{\text{th}}$ restart
\begin{eqnarray*}
     N_k = \dfrac{4\sqrt{2} d^{1/4} \an{\sqrt{M_2 M}}}{\mu R} \cdot 2^k
\end{eqnarray*}

and the number of oracle calls
\begin{eqnarray*}
     T_k = \dfrac{2048 \kappa(p, d) \lambda(p,d) d M_2^2}{\mu^2 R^2} \cdot 4^k.
\end{eqnarray*}

We obtain the total number of iterations
\begin{eqnarray*}
     N = \sum\limits_{k=1}^K N_k \leq 2^{K+1} \cdot \dfrac{4\sqrt{2} d^{1/4} \an{\sqrt{M_2 M}}}{\mu R} = \dfrac{8\sqrt{2} d^{1/4} \an{\sqrt{M_2 M}}}{\sqrt{2 \mu \varepsilon}} = O \left( \dfrac{d^{1/4} \an{\sqrt{M_2 M}}}{\sqrt{\mu \varepsilon}} \right).
\end{eqnarray*}

The total number of oracle calls (we use that $\lambda(p,d) = O(\log(d))= \tilde O(1)$):
\begin{eqnarray*}
     T = \sum\limits_{k=1}^K T_k \leq 2 \cdot 4^K \cdot \dfrac{2048 \kappa(p, d) \lambda(p,d) d M_2^2}{\mu^2 R^2} = 2 \cdot \dfrac{\mu R^2}{2\varepsilon} \cdot \dfrac{2048 \kappa(p, d) \lambda(p,d) d M_2^2}{\mu^2 R^2} = \tilde{O} \left(\dfrac{\kappa(p, d) d M_2^2}{\mu \varepsilon} \right).
\end{eqnarray*}

\subsection{Noisy Value Of Function}\label{noisy}
In this section (largely following the work \cite{dvinskikh2022gradient}) we estimate the maximum level of admissible noise $\Delta$ in general case, i.e. without simplifying assumptions about unbiasedness.

For simplicity, we consider  non-stochastic non-smooth convex optimization problem in the Euclidean proximal setup on a compact set $Q$: 
\begin{equation}\label{prblm}
    \min_{x\in Q\subseteq \mathbb{R}^d} f(x),
\end{equation}
where $f$ is  $M_2$-Lipschitz continuous.
We replace the objective by its smooth  approximation:
$
f_{\gamma}(x)\triangleq  \mathbb{E}_{u} f(x + \gamma u)
$, where $u$ is a  vector 
picked uniformly at random from the Euclidean unit ball $\{u:\|u\|_2\leq 1\}$.
From \cite{duchi2015optimal} it follows that
\begin{equation}\label{apprx}
f(x) \leq f_{\gamma}(x) \leq f(x) + \gamma M_2.
\end{equation}
Also from \cite{shamir2017optimal} (where $\Delta = 0$) we have that 
\begin{equation*}
    \nabla f_{\gamma}(x,e) = \frac{d}{2\gamma}\left(f_{\delta}(x+\gamma e) - f_{\delta}(x-\gamma e)\right)e,
\end{equation*}
where $f_{\delta} = f + \delta$ -- is the noisy value of $f$, $|\delta(x)|\le\Delta$ is a level of noise. Due to 
\cite{gasnikov2017stochastic} for all $r \in \mathbb{R}^d$
\begin{equation}\label{unbias}
 \mathbb{E}_e\langle  \left[\nabla f_{\gamma}(x, e) \right] - \nabla f_{\gamma}(x), r\rangle \lesssim \sqrt{d}{\Delta\|r\|_2}{\gamma^{-1}}
 \end{equation}
 and \cite{shamir2017optimal,beznosikov2020gradient}
 \begin{equation}\label{var_appendix}
   \mathbb{E}_e \left[ \| \nabla f_{\gamma}(x, e) - \mathbb{E}_e \nabla f_{\gamma}(x, e) \|^2_2 \right]\simeq \mathbb{E}_e \left[ \| \nabla f_{\gamma}(x, e) \|^2_2 \right] \lesssim 
   d M_2^2 + {d^2\Delta^2}{\gamma^{-2}} ,
    \end{equation}
where $e$ is random vector uniformly distributed on    the Euclidean unit sphere $\{e:\|e\|_2= 1\}$.  The r.h.s. of \eqref{unbias} is $\sqrt{d}$ better than in some other works, which used similar inequality \cite{beznosikov2020gradient,akhavan2020exploiting}. 

We say that an algorithm $\bf A$ (with $\nabla f_{\gamma}(x, e)$ oracle) is \textit{robust} for  $f_{\gamma}$ if the bias in the l.h.s. of \eqref{unbias} does not accumulate over method iterations. That is, if for $\bf A$ with $\Delta = 0$ 
$$\mathbb{E}f_{\gamma}(x^N) - \min_{x\in Q} f_{\gamma}(x)\le \Theta_A(N),$$
then with $\Delta > 0$ and (variance control)  ${d^2\Delta^2}{\gamma^{-2}} \lesssim d M^2$, see \eqref{var_appendix}:  
\begin{equation}\label{ns}
    \mathbb{E} f_{\gamma}(x^N) - \min_{x\in Q} f_{\gamma}(x)= O\left(\Theta_A(N) + \sqrt{d}{\Delta D}{\gamma^{-1}}\right),
\end{equation}
 where $D$ is a diameter of $Q$ (in \eqref{ns} we have to consider $N$ such that the first term in RHS is not smaller than the second one). 
 Many known methods are robust \cite{dvinskikh2022gradient}. In particular batched Accelerated gradient method from \cite{gorbunov2019optimal} is robust. This method was used in Theorem~\ref{FGM}.
 
Below we explain how to obtain the bound on the level of noise $\Delta$.

\textbf{Approximation.} First of all we need smoothed problem to approximate non-smooth one. For that from \eqref{apprx} we put $$\gamma = \frac{\varepsilon}{2M_2}.$$

\textbf{Variance control.} From \eqref{var_appendix} we can observe that stochastic gradient  will have the same variance (second moment) up to a numerical constant if 
$$\Delta\lesssim\frac{\gamma M_2}{\sqrt{d}}.$$

\textbf{Bias.} From \eqref{ns} we will have more restrictive condition on the level of noise
$$\Delta\lesssim\frac{\gamma\varepsilon}{D\sqrt{d}}.$$
Combination of Bias condition and Approximation condition leads to the bound
\begin{equation}\label{lon}
\Delta \lesssim\frac{\varepsilon^2}{D M_2\sqrt{d}}.
\end{equation}
The same reasoning holds for Lipschitz noise and for saddle-point problems.

It is important to note that this level of noise is maximum possible, see \cite{risteski2016algorithms} for details.

\end{document}